\documentclass[12pt,a4paper,twoside]{article}

\usepackage{graphicx}   %
\usepackage{enumerate}  %
\usepackage{latexsym}   %
\usepackage{amsmath}    %
\usepackage{amsbsy}     %
\usepackage{amsfonts}   %
\usepackage{mathrsfs}   %
\usepackage{amssymb}   %
\usepackage{array}      %
\usepackage{theorem}
\usepackage{color}
\usepackage{bm}

\usepackage{comment}

\definecolor{winered}{rgb}{0.8,0,0}

\definecolor{viola}{rgb}{0.4,0,0.7}
\definecolor{ciclamino}{rgb}{0.5,0,0.5}
\definecolor{verdescuro}{rgb}{0.0, 0.45, 0.5}
\definecolor{miogreen}{rgb}{0,0.5,0.2}
\definecolor{magenta-mio}{rgb}{1,0,1}
\definecolor{dcyan}{rgb}{0,0.4,0.8}
\definecolor{dmagenta}{rgb}{0.8,0,0.8}

\def\piegia #1{#1}
\def\ryoken #1{#1}
\def\giapie #1{#1}
\usepackage[colorlinks=true]{hyperref}
\hypersetup{urlcolor=blue, citecolor=red}

 \newtheorem{mainTh}{Main Theorem}
 \newtheorem{keyLem}{Key-Lemma}
 \newtheorem{cor}{Corollary}[section]
 \newtheorem{lem}{Lemma}[section]

 \newtheorem{defn}{Definition}[section]
 \newtheorem{rem}{Remark}[section]
 \newtheorem{notn}{Notation}
 
 \numberwithin{equation}{section}

\def\bS{{\mathbb{S}}}

\def\D{{\mathscr D}}

\def\R{{\mathbb R}}
\def\bB{{\mathbb B}}
\def\N{{\mathbb{N}}}

\def\H{{\mathcal H}}

\def\L{{\mathcal L}}

\def\ds{\displaystyle}
\def\ts{\textstyle}
\def\tm{\textmd}

\def\tm{\textmd}

\def\Sgn{\mathop{\mathrm{Sgn}}\nolimits}

\begin{document}
\thispagestyle{empty}

\pagestyle{myheadings}
\newcommand\testopari{\sc \small \piegia{Colli -- Gilardi -- Nakayashiki -- Shirakawa}}
\newcommand\testodispari{\sc\small \piegia{Quasi-linear Allen--Cahn equations with dynamic boundary conditions}}
\markboth{\testodispari}{\testopari}

\begin{center}
  {\LARGE\bf {\bf
	A class of quasi-linear \\[0.25ex] Allen--Cahn type equations \\[1.0ex] with dynamic boundary conditions\footnotemark[1]
  }}
\end{center}
\vspace{1ex}
\begin{center}
{\sc Pierluigi Colli}\footnotemark[2]\\
%Departimento di Matematica, Universit\`{a} di Pavia,\\
%via Ferrata 1, 27100, Pavia, Italy\\
e-mail: {\ttfamily pierluigi.colli@unipv.it}
\vspace{2ex}

{\sc Gianni Gilardi}\footnotemark[2]\\
%Departimento di Matematica, Universit\`{a} di Pavia,\\
%via Ferrata 1, 27100, Pavia, Italy\\
e-mail: {\ttfamily gianni.gilardi@unipv.it}
\vspace{2ex}

{\sc Ryota Nakayashiki}\footnotemark[3]\\
%Department of Mathematics and Informatics, Graduate School of Science, Chiba University\\
%1-33, Yayoi-cho, Inage-ku, Chiba, 263-8522, Japan\\
e-mail: {\ttfamily nakayashiki1108@chiba-u.jp}
\vspace{2ex}

{\sc Ken Shirakawa}\footnotemark[4]\\
%Department of Mathematics, Faculty of Education, Chiba University\\
%1-33, Yayoi-cho, Inage-ku, Chiba, 263-8522, Japan\\
e-mail: {\ttfamily sirakawa@faculty.chiba-u.jp}
\end{center}
\vspace{7ex}

\noindent
{\bf Abstract.} In this paper, we consider a class of coupled systems of PDEs, denoted by (ACE)$_{\varepsilon}$ for $ \varepsilon \geq 0 $. For each $ \varepsilon \geq 0 $, the system (ACE)$_{\varepsilon}$ consists of an Allen--Cahn type equation in a bounded spacial domain $ \Omega $, and another Allen--Cahn type equation on the smooth boundary $ \Gamma := \partial \Omega $, and besides, these coupled equations are transmitted via the dynamic boundary conditions. In particular, the equation in $ \Omega $ is derived from the non-smooth energy proposed by Visintin \piegia{in his monography ``Models of phase transitions'': hence,} the diffusion in $ \Omega $ is provided by a quasilinear form with singularity. The objective of this paper is to build a mathematical method to obtain meaningful $ L^2 $-based solutions to our systems, and to see some robustness of (ACE)$_\varepsilon$ with respect to $ \varepsilon \geq 0 $. On this basis, we will prove two Main Theorems 1 and 2, which will be concerned with the well-posedness of (ACE)$_\varepsilon$ for each $ \varepsilon \geq 0 $, and the continuous dependence of solutions to (ACE)$_\varepsilon$ for the variations of $ \varepsilon \geq 0 $, respectively.
\\[.3cm]
\piegia{{\bf Key words and phrases:} quasi-linear Allen--Cahn equation, dynamic boundary conditions,   
non-smooth energy functional, initial-boundary value problem, well-posed\-ness, continuous dependence.}

%%%%%%%%%%%%%%%%%%%%%%%%%%%%%%%%%%%%%%
\footnotetext[0]{%\hspace*{-7.4mm}
%%%%%%%%%%%%%%%%%%%%%%%%%%%%%%%%%%%%%%
\noindent $\hskip-0.3cm \empty^*$\,AMS Subject Classification 
35K55, % Nonlinear parabolic equations
35K59,  % Quasilinear parabolic equations 
82C26. % Dynamic and nonequilibrium phase transitions(general)
%%%%%%%%%%%%%%%%%%%%%%%%%%%%%%%%%%%%%%

$\empty^\dag$\,Departimento di Matematica, Universit\`{a} di Pavia, via Ferrata \piegia{5}, 27100, Pavia, Italy. \piegia{This author gratefully acknowledges some financial support from the GNAMPA (Gruppo Nazionale per l'Analisi Matematica, la Probabilit\`a e le loro Applicazioni) of INdAM (Istituto  Nazionale di Alta Matematica) and the IMATI -- C.N.R. Pavia.}

$\empty^\ddag$\,Department of Mathematics and Informatics, Graduate School of Science, Chiba University, 1-33, Yayoi-cho, Inage-ku, Chiba, 263-8522, Japan.

$\empty^\S$\,Department of Mathematics, Faculty of Education, Chiba University, 1-33, Yayoi-cho, Inage-ku, Chiba, 263-8522, Japan. This author is supported by Grant-in-Aid No. 16K05224, JSPS.
%\\
%\copyright Gakk\={o}tosho ****, GAKK\={O}TOSHO CO.,LTD
}
\newpage
\section*{Introduction}
\ \ \vspace{-3ex}

Let $ 0 < T < \infty $, $ \kappa > 0 $ and $ N \in \N $ be fixed constants. Let $ Q := (0, T) \times \Omega $ be a product set of a time-interval $ (0, T) $ and a bounded spatial domain $ \Omega \subset \R^{N} $. Let $ \Gamma := \partial \Omega $ be the boundary of $ \Omega $ with sufficient smoothness (when $ N > 1 $), and let $ {n}_{\Gamma} $ be the unit outer normal to $ \Gamma $. Besides, we put  $ \Sigma:= (0, T) \times \Gamma $.

In this paper, we fix a constant $ \varepsilon \geq 0 $ to consider the following system of PDEs, denoted by (ACE)$_\varepsilon$. 
\bigskip

\noindent 
(ACE)$_\varepsilon$:
\begin{equation}\label{1st.eq}
\partial_{t} u - \mathrm{div}\left(\frac{\nabla u}{|\nabla u|} + \kappa^2 \nabla u \right) + \beta(u) + g(u) \ni \theta\ \mbox{in $ Q $,}
\end{equation}
\begin{equation}\label{2nd.eq}
\begin{array}{c}
\ds\partial_{t} u_{\Gamma} - \varepsilon^2 {\mit \Delta}_{\Gamma} u_{\Gamma} + \bigl({\ts\frac{\nabla u}{|\nabla u|} + \kappa^2 \nabla u} \bigr)_{|_{\Gamma}} \cdot {n}_{\Gamma} + \beta_{\Gamma}(u_{\Gamma}) + g_{\Gamma}(u_{\Gamma}) \ni \theta_{\Gamma}
\\[1ex]
\mbox{and } u_{|_{\Gamma}} = u_{\Gamma} \mbox{ on $ \Sigma $,}
\end{array}
\end{equation}
\begin{equation}\label{3rd.eq}
u(0,\cdot) = u_0\ \mbox{in $ \Omega $, and}\ u_{\Gamma}(0,\cdot) = u_{\Gamma,0}\ \mbox{on $ \Gamma $.}
\end{equation}
The system (ACE)$_\varepsilon$ is a modified version of an Allen--Cahn type equation, proposed in \cite[\piegia{Chapter}~VI]{Visintin}, and the principal modifications are in the points that:
\begin{description}
\item[\tm{--}]the quasi-linear (singular) diffusion in \eqref{1st.eq} includes the regularization term $ \kappa^2 \nabla u $ with a small constant $ \kappa > 0 $;
\vspace{-1ex}
\item[\tm{--}] the boundary data $ u_{\Gamma} $ is governed by \piegia{the} dynamic boundary condition \eqref{2nd.eq}.
\end{description}
In general, ``Allen--Cahn type equation'' is a collective term to call gradient flows (systems) of governing energies, which include some double-well type potentials to reproduce the bi-stability of different phases, such as solid-liquid phases. The governing energy is called {\em free-energy,} and in the case of (ACE)$_\varepsilon$, the corresponding {free-energy} is provided as follows.
\begin{equation}\label{free-energy}
\begin{array}{ll}
\multicolumn{2}{l}{
\ds [u, u_{\Gamma}] \in H^1(\Omega) \times H^{\frac{1}{2}}(\Gamma) \mapsto \mathscr{F}_{\varepsilon}(u, u_{\Gamma})}
\\[2ex]
\qquad\qquad & := \ds \int_{\Omega} \left(|\nabla u| + \frac{\kappa^2}{2}|\nabla u|^2 + B(u) + G(u) \right) \, dx
\\[2ex]
& \quad \ds + \int_{\Gamma} \left(\frac{\varepsilon^2}{2} |\nabla_\Gamma u_{\Gamma}|^2 + B_{\Gamma}(u_{\Gamma}) + G_{\Gamma}(u_{\Gamma}) \right) \, d\Gamma \in (-\infty, \infty],
\end{array}
\end{equation}
with the effective domain:
\begin{equation*}
D(\mathscr{F}_\varepsilon) := \left\{ \begin{array}{l|l} 
[z, z_\Gamma] & \parbox{7cm}{
$ z \in H^1(\Omega) $, $ z_\Gamma \in H^{\frac{1}{2}}(\Gamma) $, $ \varepsilon z_{\Gamma} \in H^1(\Gamma) $, and $ z_{|_\Gamma} = z_\Gamma $ in $ H^{\frac{1}{2}}(\Gamma) $
}
\end{array} \right\}.
\end{equation*}
In the context, ``$ |_{\Gamma} $'' denotes the trace (boundary-value) on $ \Gamma $ for a Sobolev function, 
$ d\Gamma $ denotes the area-element on $ \Gamma $, $ \nabla_{\Gamma} $ denotes the surface gradient on $ \Gamma $, 
and  $ {\mit \Delta}_{\Gamma} $ denotes the Laplacian on the surface, i.e., the so-called Laplace-Beltrami operator. 
$ B : \R \to [0,\infty] $ and $ B_{\Gamma} : \R \to [0,\infty] $ are given proper l.s.c.\ and convex functions, 
and $ \beta = \partial B $ and $ \beta_{\Gamma} = \partial B_{\Gamma} $ are the subdifferentials of $ B $ and $ B_{\Gamma} $, respectively. 
$ G : \R \to \R $ and $ G_{\Gamma} : \R \to \R $ are $ C^1 $-functions, 
that have locally Lipschitz differentials $ g $ and $ g_{\Gamma} $, respectively. 
$ \theta : Q \to \R $ and $ \theta_{\Gamma} : \Sigma \to \R $ are given heat sources of (relative) temperature, 
and $ u_{0} : \Omega \to \R $ and $ u_{\Gamma,0} : \Gamma \to \R $ are initial data for the components $ u $ and $ u_{\Gamma} $, respectively.

In \eqref{free-energy}, the functions:
\begin{equation*}
\sigma \in \R \mapsto B(\sigma) + G(\sigma) \in (-\infty,\infty]\ \mbox{and}\ \sigma \in \R \mapsto B_{\Gamma}(\sigma) + G_{\Gamma}(\sigma) \in (-\infty,\infty],
\end{equation*}
correspond to the double-well potentials, and for instance, the setting:
\begin{equation*}
B(\sigma) = B_{\Gamma}(\sigma) = I_{[-1, 1]}(\sigma) \mbox{ and } G(\sigma) = G_{\Gamma}(\sigma) = -\frac{1}{2}\sigma^2 \mbox{, for $ \sigma \in \R $},
\end{equation*}
with use of the indicator function: 
\begin{equation*}
\sigma \in \R \mapsto I_{[-1, 1]}(\sigma) := \left\{ \begin{array}{ll}
0, & \mbox{if $ \sigma \in [-1, 1] $,}
\\[0.5ex]
\infty, & \mbox{otherwise,}
\end{array} \right.
\end{equation*}
is known as one of representative choices of the components (cf.\ \cite{Visintin}). 

Additionally, it should be noted that the presence or absence of the term 
$$ \piegia{\frac{\varepsilon^2}{2}  \int_\Gamma  |\nabla_{\Gamma} u_{\Gamma}|^2 \, d\Gamma} $$ 
brings the gap of effective domains $ D(\mathscr{F}_{\varepsilon}) $ between the cases when $ \varepsilon > 0 $ 
and $ \varepsilon = 0 $. More precisely, the domains $ D(\mathscr{F}_{\varepsilon}) $ when $ \varepsilon >0 $ will 
uniformly \piegia{coincide} with a convex subset in $ H^1(\Omega) \times H^1(\Gamma) $, and this convex set will be a 
proper subset of the domain $ D(\mathscr{F}_0) $ when $ \varepsilon = 0 $ which will be located in the 
\piegia{wider space} $ H^1(\Omega) \times H^{\frac{1}{2}}(\Gamma) $. 

In the case when the diffusion in \eqref{1st.eq} is just given by the usual Laplacian, the corresponding Allen--Cahn equation has been studied by a number of researches \piegia{(cf., e.g., \cite{CC13,CF1,GG08,GGM08,Is})}, 
%from various viewpoints. In line of the 
and some qualitative results for $ L^2 $-based solutions were obtained 
by means of the theories of parabolic PDEs, in \cite{LSU,LM}. \piegia{To investigate dynamic boundary conditions,
our approach exploits techniques similar to those employed in \cite{CC13} and resumed in other 
solvability studies and optimal control theories, the reader may see 
\cite{CF1,CF2,CF3,CGS0,CGS1,CGS2,CS}. 
Still about dynamic boundary conditions, let us point out that there has been a recent growing interest about the justification and the study of phase field models, as well as systems of Allen--Cahn and Cahn--Hilliard type, 
including dynamic boundary conditions. Without trying to be exhaustive, let us 
mention at least the papers \cite{ChGaMi, GG08, GaWa, GiMiSc, GMS10,  GoMi, GMS11, Li, MRSS}.}

Nevertheless, the mathematical analysis for our system (ACE)$_{\varepsilon}$ will not be just an analogy work with the previous ones. 
In fact, due to the singularity of the diffusion $ -\mathrm{div}(\frac{\nabla u}{|\nabla u|} + \kappa^2 \nabla u) $ in \eqref{1st.eq}, it will not be so easy to apply the theories of \cite{LSU,LM}, and to see the $ L^2 $-based expression of the first variation of the {free-energy}.

In view of this, we set the goal in this paper to show the following two Main Theorems, which are concerned with qualitative properties of the systems (ACE)$_{\varepsilon}$ for $ \varepsilon \geq 0 $.
\begin{description}
\item[Main Theorem 1:] the well-posedness for (ACE)$_{\varepsilon}$, for all $ \varepsilon \geq 0 $.
\vspace{-1ex}
\item[Main Theorem 2:] the continuous dependence of solutions to (ACE)$_{\varepsilon}$ with respect to the value of $ \varepsilon \geq 0 $, and especially the (\piegia{right-hand}) continuity at $ \varepsilon = 0 $.
\end{description}

The content of this paper is as follows. 
The Main Theorems 1 and 2 are stated in Section 2, and these results are discussed on the basis of the preliminaries prepared in Section 1, and Key-Lemmas 1--3. Based on this, we give the proofs of the Key-Lemmas and Main Theorems in the remaining Sections 4 and 5, respectively.
 
\section{Preliminaries}
\ \ \vspace{-3ex}

In this Section, we outline some basic notations and known facts, as preliminaries of our study. 

\begin{notn}[Notations in real analysis]\label{Note00}
\begin{em}
For arbitrary $ a, b \in [-\infty, \infty] $, we define:
\begin{equation*}
a \vee b := \max \{ a, b \} \mbox{ and } a \wedge b := \min \{ a, b \}.
\end{equation*}

Let $ d \in \N $ be any fixed dimension. Then, we simply denote by $ |x| $ and $ x \cdot y $ the Euclidean norm of $ x \in \mathbb{R}^d $ and the standard scalar product of  $ x, y \in \R^d $, \piegia{respectively.}
%, i.e.:
%\begin{equation*}
%\begin{array}{c}
%\displaystyle
%| x | := {\textstyle \sqrt{x_1^2 + \cdots +x_d^2}} \mbox{ and } x \cdot y  := x_1 y_1 + \cdots +x_d y_d, 
%\\[1ex]
%\mbox{ for all $ x = [x_1, \dots, x_d], \ y = [y_1, \dots, y_d] \in \mathbb{R}^d $.}
%\end{array}
%\end{equation*}
Besides, we denote by $ \bB^d $ and $ \bS^{d -1} $ the $ d $-dimensional unit open ball centered at the origin, 
and its boundary, \piegia{respectively.}
%, i.e.:
%\begin{equation*}
%\bB^d := \left\{ \begin{array}{l|l}
%x \in \R^d & |x| < 1
%\end{array} \right\} \mbox{ and } \bS^{d -1} := \left\{ \begin{array}{l|l}
%x \in \R^d & |x| = 1
%\end{array} \right\}.
%\end{equation*}

For any $ d \in \N $, the $ d $-dimensional Lebesgue measure is denoted by $ \L^d $, and $ d $-di\-men\-sion\-al Hausdorff meausure is denoted by~$ \H^d $. 
Unless otherwise specified, the measure theoretical phrases, 
such as ``a.e.'', ``$ dt $'', ``$ dx $'', and so on, are  with respect to the Lebesgue measure in each corresponding dimension. 
Also, in the observations on a smooth \piegia{surface~$ S $}, the phrase ``a.e.'' 
is with respect to the Hausdorff measure in each corresponding Hausdorff dimension, 
and the area element on $ S $ is denoted by~$ dS $. 

Additionally, we mention about the following elementary fact, which is used, frequently, in the proofs of Key-Lemmas and Main Theorems.
\begin{description}
\item[(Fact\,0)] Let $ m \in \N $ be a fixed finite number. If $ \{ \alpha_1, \dots, \alpha_m \} \subset \R $ and $ \{ a_n^k \}_{n = 1}^{\infty} $, $ k= 1,\dots, m $, fulfill \piegia{that}
\begin{equation*}
\varliminf_{n \to \infty} a_n^k \geq \alpha_k \mbox{, $ k=1, \dots, m $, and } \varlimsup_{n \to \infty} \sum_{k=1}^{m} a_n^k \leq \sum_{k=1}^m \alpha_k
\end{equation*}
\piegia{then}, it holds \piegia{that}
\begin{equation*}
\lim_{n \to \infty} a_n^k = \alpha_k \mbox{, $ k=1, \dots, m $.}
\end{equation*}
\end{description}

\begin{comment}
For a (Lebesgue) measurable function $ f : B \rightarrow [-\infty, \infty] $ on a Borel subset $ B \subset \R^d $, we denote by $ [f]^+ $ and $ [f]^- $, respectively, the positive part and the negative part of $ f $, i.e.:
\begin{equation*}
[f]^+(x) := f(x) \vee 0 \mbox{ \ and \ } [f]^-(x) := - (f(x) \wedge 0), \mbox{ for a.e. $ x \in B $.}
\end{equation*}
Let $ U \subset \mathbb{R}^N $ be an open set. We denote by $ \mathcal{M}(U) $ the space of all finite Radon measures on $ U $. We recall that the space $ \mathcal{M}(U) $ is  the dual space of the Banach space $ C_0(U) $, where $ C_0(U) $ denotes the closure of the space $ C_{\rm c}(U) $ of all continuous functions having compact supports, in the topology of $ C(\overline{U}) $.
\end{comment}
\end{em}
\end{notn}

\begin{notn}[Notations of functional analysis]\label{NoteFuncSp}
\begin{em}
For an abstract Banach space $ X $, we denote by $ |{}\cdot{}|_X $ the norm of $ X $, and denote by $ {}_{X^*} \langle{}\cdot{}, {}\cdot{}\rangle_X $ the duality pairing between $ X $ and the dual space $ X^* $ of $ X $. Let $ \mathcal{I}_X : X \to X $ be the identity map from $ X $ onto $ X $. In particular, when $ X $ is a Hilbert space, we denote by $ ({}\cdot{},{}\cdot{})_X $ the inner product in $ X $. 

\begin{comment}
We denote by $ F : H^1(\Omega) \rightarrow H^1(\Omega)^* $ the duality map, defined as:
\begin{equation*}
\langle Fz, w \rangle_{H^1(\Omega)^*, H^1(\Omega)} := \int_\Omega \nabla z \cdot \nabla w \, dx +\int_\Omega zw \, dx, 
\mbox{ \ for $ [z, w] \in H^1(\Omega)^2 $.}
\end{equation*}
\end{comment}
\piegia{Here and in the sequel, $\Omega$ denotes an open subset of~$\R^N$, which we assume to be bounded and smooth.
Moreover, $ \Gamma $ and $n_\Gamma$ denote its boundary~$ \partial \Omega $ and the outward unit normal vector field on~$\Gamma$, respectively.}
Let $ {\mit \Delta}_{N} $ be the \piegia{Laplace operator}, subject to 
the Neumann-zero boundary condition, which is defined as:
\begin{equation}
\label{gianni1}
\begin{array}{rcl}
{\mit \Delta}_{N} : v \in D({\mit \Delta}_N) & := & \left\{ \begin{array}{l|l}
z \in H^2(\Omega) & (\nabla z)_{|_\Gamma} \cdot \piegia{n_{\Gamma}} = 0 \mbox{ in $ H^{\frac{1}{2}}(\Gamma) $}
\end{array} \right\} \subset L^2(\Omega) 
\\[1ex]
& \mapsto & {\mit \Delta}_N v := {\mit \Delta} v \in L^2(\Omega).
\end{array}
\end{equation}
In this paper, we identify the \piegia{unbounded closed} operator $ -{\mit \Delta}_N $ \piegia{with its} linear and continuous \piegia{extension} from $ H^1(\Omega) $ into $ H^{1}(\Omega)^* $, \piegia{by setting:}
\begin{equation*}
\piegia{{}_{ H^{1}(\Omega)^*} \langle{}  {\mit \Delta}_N z        {}, {}     w      {}\rangle_{ H^{1}(\Omega)} }
= \int_\Omega \nabla z \cdot \nabla w \, dx, ~ \mbox{ for all $ [z, w] \in \piegia{H^{1}(\Omega) \times H^{1}(\Omega)} $.}
\end{equation*}
\end{em}
\end{notn}

\begin{rem}\label{Rem.zero}
\begin{em}
Note that the boundary $ \Gamma = \partial \Omega $ has no boundary as $ (N -1) $-dimensional surface. 
Therefore, for any $ s > 0 $, \piegia{the dual space $ H^{-s}(\Gamma) := H^s(\Gamma)^* $ of the Sobolev space $ H^s(\Gamma) $ 
coincides with the closure of the class of the smooth functions on $\Gamma$ in the topology of~$ H^s(\Gamma)^* $}.
\end{em}
\end{rem}

\begin{notn}[Notations of surface-differentials]\label{bdryOp}
\begin{em}
\begin{comment}
Let $ [\,\cdot,n_\Gamma]_{_\Gamma} $ be a linear operator from $ C^\infty(\overline{\Omega}) $ into $ \mathscr{D}'(\Gamma) $, which is defined as:
\begin{equation}\label{sNormal}
\begin{array}{c}
\ds \langle [\varpi, n_\Gamma]_{_\Gamma}, \varphi \rangle := \int_\Omega {\rm div} \, \varpi \, \varphi^{\rm ex} \, dx + \int_\Omega \varpi \cdot \nabla \varphi^{\rm ex} \, dx,
\\[2ex]
\mbox{for all $ \varpi \in C^\infty(\overline{\Omega})^N $ and $ \varphi \in C^\infty(\Gamma) $,}
\end{array}
\end{equation}
by using the extension $ \varphi^{\rm ex} \in C^\infty(\overline{\Omega}) $ of each $ \varphi \in C^\infty(\Gamma) $.
\end{comment}
\piegia{Since $\Omega$ is bounded and smooth, there exists a function $d_\Gamma:\overline\Omega\to\R$ such~that\piegia{%
%\begin{enumerate}
%\item[$(*)$]
\begin{equation}
\label{gianni2}
\hbox{$ d_\Gamma$ is smooth on $\overline\Omega$ \ and \ 
$d_\Gamma(x) = \displaystyle \inf_{y \in \Gamma} |x -y| $ for $x$ in a neighborhood of~$\Gamma$.}
\end{equation}
%\end{enumerate}
}%
We notice that $n_\Gamma(x)=-\nabla d_\Gamma(x)$ for every $x\in\Gamma$.}
On this basis, let $ \nabla_\Gamma $ be the operator of surface-gradient on $ \Gamma $, which is defined as:
\begin{equation}\label{sGrad}
\nabla_\Gamma : \varphi \in \piegia{{}C^1{}}(\Gamma) \mapsto
  \nabla_\Gamma \varphi := \nabla \varphi^{\rm ex} -(\nabla {d}_\Gamma \otimes \nabla {d}_\Gamma) \nabla \varphi^{\rm ex}
  \in \piegia{{}C^0{}}(\Gamma)^N,
\end{equation}
by using the extension $ \varphi^{\rm ex} \in \piegia{{}C^1{}}(\overline{\Omega}) $ of each $ \varphi \in \piegia{{}C^1{}}(\Gamma) $. 
Let $ {\rm div}_\Gamma $ be the operator of surface-divergence, which is defined as:
\begin{equation}\label{sDiv}
{\rm div}_\Gamma : \omega \in \piegia{{}C^1{}}(\Gamma)^N \mapsto {\rm div}_\Gamma \omega := {\rm div} \omega^{\rm ex} - \nabla (\omega^{\rm ex} \cdot \nabla d_{\Gamma}) \cdot \nabla d_{\Gamma} \in \piegia{{}C^0{}}(\Gamma),
\end{equation}
by using the extension $ \omega^{\rm ex} \in \piegia{{}C^1{}}(\overline{\Omega})^N $ of each $ \omega \in \piegia{{}C^1{}}(\Gamma)^N $.

It is known that the definition formulas \eqref{sGrad}--\eqref{sDiv} are well-defined, 
and the values $ \nabla_\Gamma \varphi $ and $ {\rm div}_\Gamma \omega $ are settled 
independently of the choices of \piegia{the} extensions 
$ \varphi^{\rm ex} \in \piegia{{}C^1{}}(\overline{\Omega}) $ and $ \omega^{\rm ex} \in \piegia{{}C^1{}}(\overline{\Omega})^N $ 
of $ \varphi \in \piegia{{}C^1{}}(\Gamma) $ and $ \omega \in \piegia{{}C^1{}}(\Gamma)^N $, respectively,
\piegia{and of the function $d_\Gamma$ \piegia{satisfying~\eqref{gianni2}}}. 

On the basis of \eqref{sGrad}--\eqref{sDiv}, the Laplace--Beltrami operator $ {\mit \Delta}_\Gamma $, i.e., 
the surface-Laplacian on $ \Gamma $ is defined as follows:
\begin{equation*}
{\mit \Delta}_\Gamma : \varphi \in \piegia{{}C^2{}}(\Gamma) \mapsto {\mit \Delta}_\Gamma \varphi := {\rm div}_\Gamma (\nabla_\Gamma \varphi) \in \piegia{{}C^0{}}(\Gamma),
\end{equation*}
\piegia{by using the extension $ \varphi^{\rm ex} \in \piegia{{}C^2{}}(\overline{\Omega}) $ of each $ \varphi \in \piegia{{}C^2{}}(\Gamma) $.} 
\end{em}
\end{notn}

\begin{rem}\label{Rem.sOps}
\begin{em}
Let us define a closed subspace $ \bm{L}_{\rm div}^2(\Omega) $ in $ L^2(\Omega)^N $ and a closed subspace $ \bm{L}_{\rm tan}^2(\Gamma) $ in $ L^2(\Gamma)^N $, by putting:
\begin{equation*}
\bm{L}_{\rm div}^2(\Omega) := \left\{ \begin{array}{r|l}
\omega \in L^2(\Omega)^N & {\rm div} \, \omega \in L^2(\Omega)
\end{array} \right\}, \mbox{ and}
\end{equation*}
\begin{equation*}
\bm{L}_{\rm tan}^2(\Gamma) := \left\{ \begin{array}{r|l}
\omega \in L^2(\Gamma)^N & \omega \cdot n_\Gamma = 0 \mbox{ a.e. on $ \Gamma $}
\end{array} \right\}, \mbox{ respectively.}
\end{equation*}
Then, on account of the general theories as in \cite{Kenmochi75,SV97}, we can see the following 
facts \piegia{(cf.\ \cite{Kenmochi75})}.
\begin{description}
\item[(Fact\,1)] The mapping \ $ 
\nu \in \piegia{H^{1}(\Omega)^N} \mapsto \nu_{|_\Gamma} \cdot n_\Gamma \in H^{\frac{1}{2}}(\Gamma)
$ \ can be extended as a linear and continuous operator $ [\,(\cdot) \cdot n_\Gamma]_{_\Gamma} $ from $ \bm{L}_{\rm div}^2(\Omega) $ into $ H^{-\frac{1}{2}}(\Gamma) $, such that:
\begin{equation*}
\begin{array}{c}
\ds {}_{{}^{H^{-\piegia{1/2}}(\Gamma)}}\langle [\nu \cdot n_\Gamma]_{_\Gamma}, z_{|_\Gamma} \rangle_{{}^{H^{\piegia{1/2}}(\Gamma)}}
  = \int_\Omega {\rm div}\,\nu \, z \, dx +\int_\Omega \nu \cdot \nabla z \, dx,
\\[2ex]
\mbox{for all $ \nu \in \bm{L}_{\rm div}^2(\Omega) $ and $ z \in H^1(\Omega) $.}
\end{array}
\end{equation*}
\item[(Fact\,2)]The surface gradient $ \nabla_\Gamma $ can be extended as a linear and continuous operator from $ H^1(\Gamma) $ into $ \bm{L}_{\rm tan}^2(\Gamma) $. The extension is derived in the definition process of the space $ H^1(\Gamma) $ as the completion of $ \piegia{{}C^1{}}(\Gamma) $. Then, the topology of the completion is taken with respect to the norm, induced by the following bi-linear form:
\begin{equation*}
[\varphi, \psi] \in \piegia{{}C^1{}}(\Gamma)^2 \mapsto \int_\Gamma \bigl( \varphi \psi +\nabla_\Gamma \varphi \cdot \nabla_\Gamma \psi \bigr) \, d\Gamma.
\end{equation*}
The inner product in $ (\cdot,\cdot)_{H^1(\Gamma)} $ is given as the extension of the above bi-linear form. Hence, in this paper, we identify the operator $ \nabla_\Gamma $ with the extension from $ H^{1}(\Gamma) $ into $ \bm{L}_{\rm tan}^2(\Gamma) $. 

%Additionally, since the boundary $ \Gamma = \partial \Omega $ has no boundary as $ (N -1) $-dimensional surface, the dual space $ H^1(\Gamma)^* $ of $ H^1(\Gamma) $ coincides with the closure $ H^{-1}(\Gamma) $ $ (= \overline{C_{\rm c}^\infty(\Omega)}) $ of the class $ C_{\rm c}^\infty(\Gamma) $ of smooth functions with compact supports on $ \Gamma $ in the topology of $ H^1(\Gamma)^* $. 
\item[(Fact\,3)]The surface divergence $ {\rm div}_\Gamma $ can be extended as a linear and continuous operator from $ \bm{L}_{\rm tan}^2(\Gamma) $ into $ H^{-1}(\Gamma) $ $ (= H^1(\Gamma)^*) $, via the following Green-type formula:
\begin{equation*}
\begin{array}{c}
\ds\piegia{{}_{ H^{-1}(\Gamma)} \langle   \,  {\rm div}_\Gamma w ,  z \,     \rangle_{ H^1(\Gamma)}   =  - {}}\int_\Gamma w \cdot \nabla_\Gamma z \, d \Gamma \mbox{, \ for \piegia{all $ z \in H^1(\Gamma) $ and $ w \in \bm{L}_{\rm tan}^2(\Gamma) $.}}
%\\[2ex]
%\mbox{any $ w \in \bm{L}_{\rm tan}^2(\Gamma) $ satisfying $ {\rm div}_\Gamma w \in L^2(\Gamma) $.}
\end{array}
\end{equation*}
Hence, in this paper, we regard the Laplace--Beltrami operator $ {\mit \Delta}_\Gamma = {\rm div}_\Gamma \circ \nabla_\Gamma $ as \piegia{a} linear and continuous operator from $ H^1(\Gamma) $ into $ H^{-1}(\Gamma) $. In particular, the operator $ -{\mit \Delta}_\Gamma $ forms a duality map between $ H^1(\Gamma) $ and $ H^{-1}(\Gamma) $. 
\end{description}
\end{em}
\end{rem}

\begin{notn}[Notations in convex analysis]\label{NoteConvex}
\begin{em}
For any proper lower semi-continuous (l.s.c.\ from now on) and convex function $ \Psi $ defined on a Hilbert space $ X $, we denote by $ D(\Psi) $ its effective domain, and denote by $ \partial \Psi $ its subdifferential. The subdifferential $ \partial \Psi $ is a set-valued map corresponding to a weak differential of $ \Psi $, and it \piegia{turns out to be} a maximal monotone graph in the product space $ X \times X $. More precisely, for each $ z_0 \in X $, the value $ \partial \Psi(z_0) $ is defined as a set of all elements $ z_0^* \in X $ which satisfy the following variational inequality:
\begin{equation*}
(z_0^*, z -z_0)_X \leq \Psi(z) -\Psi(z_0) \mbox{, for any $ z \in D(\Psi) $.}
\end{equation*}
The set \
$ D(\partial \Psi) := \{ z \in X \,|\, \partial \Psi(z) \ne \emptyset \} $ \
is called the domain of $ \partial \Psi $. We often use the notation ``$ [z_0, z_0^*] \in \partial \Psi $ in $ X \times X $'', to mean that ``$ z_0^* \in \partial \Psi(z_0) $ in $ X $ with $ z_0 \in D(\partial \Psi) $'', by identifying the operator $ \partial \Psi $ with its graph in $ X \times X $. \piegia{Let us refer to \cite{Attouch,Barbu,Brezis, Kenmochi81} for definitions, properties, results about subdifferentials and maximal monotone operators}. 
\end{em}
\end{notn}

\begin{rem}\label{exConvex}
\begin{em}
As one of representatives of the subdifferentials, we exemplify the following set-valued function $ {\rm Sgn} : \R^N \rightarrow 2^{\R^N} $, given as:
\begin{equation*}
\omega \in \R^N \mapsto {\rm Sgn}(\omega) := \left\{ \begin{array}{ll}
\ds \frac{\omega}{|\omega|}, & \mbox{if $ \omega \ne 0 $,}
\\[2.5ex]
\overline{\bB^N}, & \mbox{otherwise.}
\end{array} \right.
\end{equation*}
It is known that the set-valued function $ {\rm Sgn} $ coincides with the subdifferential of the Euclidean norm \piegia{$ |{}\cdot{}| : \omega \in \R^N \mapsto |\omega| = \sqrt{\omega \cdot \omega} \in [0, \infty) $.}
%, i.e.:
%\begin{equation*}
%\partial |{}\cdot{}|(\omega) = {\rm Sgn}(\omega), \mbox{ for any $ \omega \in D(\partial |{}\cdot{}|) = \R^N $.}
%\end{equation*}
Also, it is known that (cf.\ \cite{Barbu,Brezis}) the operator $ -{\mit \Delta}_N $ 
\piegia{defined in \eqref{gianni1}}
coincides with the subdifferential of \piegia{the}  proper l.s.c.\ and convex function $ \Psi_N $ on $ L^2(\Omega) $, defined as:
\begin{equation*}
z \in L^2(\Omega) \mapsto \Psi_N(z) := \left\{ \begin{array}{ll}
\multicolumn{2}{l}{\ds \frac{1}{2} \int_\Omega |\nabla z|^2 \, dx, \mbox{ if $ z \in H^1(\Omega) $,}}
\\[2ex]
\infty, & \mbox{otherwise.}
\end{array} \right.
\end{equation*}
More precisely\piegia{, we have}:
\begin{equation*}
\partial \Psi_N(z) = \{ -{\mit\Delta}_{N} z \} \mbox{ in $ L^2(\Omega) $, for any $ z \in D(\partial \Psi_N) = D({\mit \Delta}_N) $.}
\end{equation*}
\end{em}
\end{rem}

\begin{rem}[Time-dependent subdifferentials]\label{Rem.Time-dep.}
\begin{em}
It is often useful to consider the subdifferentials under time-dependent settings of convex functions. 
With regard to this topic, certain general theories were established by \piegia{some} researchers (e.g., Kenmochi~\cite{Kenmochi81} and \^{O}tani~\cite{Otani}). So, referring to\piegia{, e.g., \cite[Chapter~2]{Kenmochi81} or \cite[Remark 1.1 (Fact\,1)]{SWY14},} we can see the following fact. 
\begin{description}
\item[{(\hypertarget{Fact3}{Fact\,4})}]Let $ E_0 $ be a convex subset in a Hilbert space $ X $, let $ I \subset [0, \infty) $ be a time-interval, and for any $ t \in I $, let $ \Psi^t : X \rightarrow (-\infty, \infty] $ be a proper l.s.c.\ and convex function such that $ D(\Psi^t) = E_0 $ for all $ t \in I $. Based on this, let us define a convex function $ \hat{\Psi}^I : L^2(I; X) \rightarrow (-\infty, \infty] $, by putting:
\begin{equation*}
\zeta \in L^2(I; X) \mapsto \hat{\Psi}^I(\zeta) := \left\{ \begin{array}{ll}
\multicolumn{2}{l}{\ds \int_I \Psi^t(\zeta(t)) \, dt, \mbox{ if $ \Psi^{(\cdot)}(\zeta) \in L^1(I) $,}}
\\[1ex]
\infty, & \mbox{otherwise.}
\end{array} \right.
\end{equation*}
Here, if $ E_0 \subset D(\hat{\Psi}^I) $, i.e., if the function $ t \in I \mapsto \Psi^t(z) $ is integrable for any $ z \in E_0 $, then it holds that:
\begin{equation*}
\begin{array}{c}
[\zeta, \zeta^*] \in \partial \hat{\Psi}^I \mbox{ in $ L^2(I; X) \times L^2(I;X) $, \piegia{iff}}
\\[1ex]
\zeta \in D(\hat{\Psi}^I) \mbox{ and } [\zeta(t), \zeta^*(t)] \in \partial \Psi^t \mbox{ in $ X \times X $, a.e. $ t \in I $.}
\end{array}
\end{equation*}
\end{description}
\end{em}
\end{rem}

Finally, we mention about notions of  \piegia{convergence for functionals}.
 
\begin{defn}[Mosco convergence: cf.\ \cite{Mosco}]\label{Def.Mosco}
\begin{em}
Let $ X $ be an abstract Hilbert space. Let $ \Psi : X \rightarrow (-\infty, \infty] $ be a proper l.s.c.\ and convex function, and let $ \{ \Psi_n \}_{n = 1}^\infty $ be a sequence of proper l.s.c.\ and convex functions $ \Psi_n : X \rightarrow (-\infty, \infty] $, $ n \in \N $.  Then, it is said that $ \Psi_n \to \Psi $ on $ X $, in the sense of Mosco~\cite{Mosco}, as $ n \to \infty $, \piegia{iff} the following two conditions are fulfilled.
\begin{description}
\item[(\hypertarget{M_lb}{M1}) The condition of lower-bound:]$ \ds \varliminf_{n \to \infty} \Psi_n(z_n^\dag) \geq \Psi(z^\dag) $, if $ z^\dag \in X $, $ \{ z_n^\dag  \}_{n = 1}^\infty \subset X $, and $ z_n^\dag \to z^\dag $ weakly in $ X $ as $ n \to \infty $.
\item[(\hypertarget{M_opt}{M2}) The condition of optimality:]for any $ z^\ddag \in D(\Psi) $, there exists a sequence \linebreak $ \{ z_n^\ddag \}_{n = 1}^\infty \subset X $ such that $ z_n^\ddag \to z^\ddag $ in $ X $ and $ \Psi_n(z_n^\ddag) \to \Psi(z^\ddag) $, as $ n \to \infty $.
\end{description} 
\end{em}
\end{defn}

\begin{rem}\label{MY}
\begin{em}
(cf.\ \cite[Proposition 2.68 and Theorem 3.26]{Attouch}) For a proper l.s.c.\ and convex function $ \Psi : H \to (-\infty,\infty] $ on a Hilbert space $ H $, it is known that the sequence $ \{ \Psi^{\lambda} \}_{\lambda > 0} $ of Moreau--Yosida regularizations:
\begin{equation*}
z \in H \mapsto \Psi^{\lambda}(z) := \inf \left\{ \begin{array}{l|l} \ds\frac{1}{2\lambda} | \tilde{z} - z |_{H}^2 + \Psi(\tilde{z}) & \tilde{z} \in H \end{array}\right\} \mbox{, for $ \lambda > 0 $,}
\end{equation*}
converges to $ \Psi $ on $ H $, in the sense of \piegia{Mosco}, as $ \lambda \downarrow 0 $.
\end{em}
\end{rem}

\begin{defn}[$\Gamma$-convergence: cf.\ \cite{GammaConv}]\label{Def.Gamma}
\begin{em}
Let $ X $ be an abstract Hilbert space, $ \Psi : X \to (-\infty, \infty] $ be a proper functional, and $ \{ \Psi_n \}_{n = 1}^\infty $ be a sequence of proper functionals $ \Psi_n : X \rightarrow (-\infty, \infty] $, $ n \in \N $. We say that $ \Psi_n \to \Psi $ on $ X $, in the sense of $ \Gamma $-convergence~\cite{GammaConv}, as $ n \to \infty $ \piegia{iff} the following two conditions are fulfilled.
\begin{description}
\item[(\hypertarget{G_lb}{\boldmath $\Gamma$1}) The condition of lower-bound:]$ \ds \varliminf_{n \to \infty} \Psi_n(z_n^\dag) \geq \Psi(z^\dag) $ if $ z^\dag \in X $, $ \{ z_n^\dag \}_{n = 1}^\infty \subset X $, and $ z_n^\dag \to z^\dag $ (strongly) in $ X $ as $ n \to \infty $.
\item[(\hypertarget{G_opt}{\boldmath $\Gamma$2}) The condition of optimality:]for any $ z^\ddag \in D(\Psi) $, there exists a sequence\linebreak  $ \{ z_n^\ddag  \}_{n = 1}^\infty \subset X $ such that $ z_n^\ddag \to z^\ddag $ in $ X $ and $ \Psi_n(z_n^\ddag) \to \Psi(z^\ddag) $ as $ n \to \infty $.
\end{description} 
\end{em}
\end{defn}

\begin{rem}\label{Rem.MG}
\begin{em}
\piegia{Of course, the $ \Gamma $-convergence recalled in Definition~\ref{Def.Gamma} 
is the one associated with the strong topology of $X$.}
Note that if the functionals are convex, then \piegia{the Mosco convergence just introduced} 
implies $ \Gamma $-convergence, i.e., the $ \Gamma $-convergence of convex functions can be regarded as a weak version of Mosco convergence. 
Additionally, as a basic matter of the \piegia{Mosco convergence}, we can see the following fact 
(see \cite[Chapter~2]{Kenmochi81}, or \cite[Remark 1.5 (Fact\,7)]{SWY14}, for example).
\begin{description}
\item[(\hypertarget{Fact4}{Fact\,5})]Let $ X $, $ \Psi $ and $ \{ \Psi_n \}_{n = 1}^\infty $ \piegia{be as in 
Definition~\ref{Def.Mosco}%{\bf (INSTEAD OF DEFINITION~1.2 SINCE WE SUSPECT THAT HERE WE NEED CONVEXITY OF FUNCTIONS)}
}. Besides, let us assume that:
\begin{center}
$ \Psi_n \to \Psi $ on $ X $, in the sense of  \piegia{$ \Gamma $-convergence}, as $ n \to \infty $,
\vspace{-1ex}
\end{center}
and
\begin{equation*}
\left\{ ~ \parbox{9cm}{
$ [z, z^*] \in X \times X $, ~ $ [z_n, z_n^*] \in \partial \Psi_n $ in $ X \times X $, $ n \in \N $,
\\[1ex]
$ z_n \to z $ in $ X $ and $ z_n^* \to z^* $ weakly in $ X $, as $ n \to \infty $.
} \right.
\end{equation*}
Then, it holds that:
\begin{equation*}
[z, z^*] \in \partial \Psi \mbox{ in $ X \times X $, and } \Psi_n(z_n) \to \Psi(z) \mbox{, as $ n \to \infty $.}
\end{equation*}
\end{description}
\end{em}
\end{rem}

\section{Statements of Main Theorems}
\ \ \vspace{-3ex}

First, we configure the base-space of solutions to the systems (ACE)$_\varepsilon$, for $ \varepsilon \geq 0 $. In any case of $ \varepsilon \geq 0 $, the base-space is settled by a product Hilbert space: 
\begin{equation*}
\mathscr{H} := L^2(\Omega) \times L^2(\Gamma),
\end{equation*}
endowed with the inner product:
\begin{equation*}
\begin{array}{c}
\left( [z_1, z_{\Gamma, 1}], [z_2, z_{\Gamma, 2}] \right)_{\mathscr{H}} := (z_1, z_2)_{L^2(\Omega)} +(z_{\Gamma, 1}, z_{\Gamma, 2})_{L^2(\Gamma)}, 
\\[2ex]
\mbox{for any $ [z_k, z_{\Gamma, k}] $, $ k = 1, 2 $.}
\end{array}
\end{equation*}

Next, we prescribe the assumptions in our study.
\begin{description}
\item[(A0)]$ \piegia{N \in \N} $ and $ 0 < T < \infty $ are fixed constants, 
and $ \Omega $ is a bounded domain \piegia{in $\R^N$} with a smooth boundary~$ \Gamma $.
\piegia{In particular, it} fulfills the \piegia{condition~\eqref{gianni2}} in Notation~\ref{bdryOp}.
\item[(A1)] $ B : D(B) \subset \R \to [0,\infty] $ and $ B_{\Gamma} : D(B_{\Gamma}) \subset \R \to [0,\infty] $ are proper l.s.c.\ and convex functions and $ \beta = \partial B \subset \R \times \R $ and $ \beta_{\Gamma} = \partial B_{\Gamma} \subset \R \times \R $ are the subdifferentials of $ B $ and $ B_{\Gamma} $, respectively. Furthermore, the convex functions $ B $ and $ B_{\Gamma} $, and the subdifferentials $ \beta $ and $ \beta_{\Gamma} $ fulfill the following conditions:
\begin{description}
\item[(a1)] $ B(0) = 0 $ and $ B_{\Gamma}(0) = 0 $, and hence $ [0,0] \in \beta $ and $ [0,0] \in \beta_{\Gamma} $ on $ \R^2 $;

\item[(a2)] there exists an interval $ I_{B} \subset \R $, such that:
\begin{equation*}
\mathrm{int} I_{B} \neq \emptyset,\ D(\beta) = D(\beta_{\Gamma}) = I_{B} \mbox{, and } B, B_\Gamma \in C(\overline{I_{B}}) \cap L^{\infty}(I_B);
\end{equation*}
\item[(a3)] there \piegia{exist} positive constants $ a_{k}, b_{k} $, $ k = 0,1 $, such that:
\begin{equation*}
a_0 \left|[\beta_{\Gamma}]^{\circ}(\tau) \right| - b_0 \leq \left|[\beta]^{\circ}(\tau)\right| \leq a_1 \left|[\beta_{\Gamma}]^{\circ}(\tau) \right| + b_1 \mbox{, for any $ \tau \in I_{B} $,}
\end{equation*}
where $ [\beta]^{\circ} $ and $ [\beta_{\Gamma}]^{\circ} $ are the minimal sections for $ \beta $ and $ \beta_{\Gamma} $, respectively.
\end{description}

\item[(A2)] $ G : \R \to \R $ and $ G_{\Gamma} : \R \to \R $ are $ W_{\mathrm{loc}}^{2,\infty} $-functions such that the differentials $ g = G' $ and $ g_{\Gamma} = \piegia{{G_{\Gamma}}'} $ are Lipschitz continuous on $ \overline{I_{B}} $.

\item[(A3)] The forcing pair $ [\theta, \theta_{\Gamma}] $ belongs to $ L^2(0,T;\mathscr{H}) $ and the initial pair $ [u_0, u_{\Gamma,0}] $ belongs to a class $ \D_{*} $, defined as:
\begin{equation}\label{D0}
\D_{*} := \left\{\begin{array}{l|l} [z,z_{\Gamma}] \in \mathscr{H} & \mbox{ $ z \in \overline{I_{B}} $ a.e. in $ \Omega $ and $  z_\Gamma \in \overline{I_B} $ a.e. on $ \Gamma $ }
\end{array} \right\}.
\end{equation}
\end{description}

In addition, let us set:
\begin{equation*}
\mathscr{V}_\varepsilon := \left\{ \begin{array}{l|l}
[z, z_\Gamma] \in \mathscr{H} & \parbox{7cm}{
$ z \in H^1(\Omega) $, $ z_\Gamma \in H^{\frac{1}{2}}(\Gamma) $, $ \varepsilon z_\Gamma \in H^1(\Gamma) $, and $ z_{|_{\Gamma}} = z_\Gamma $ a.e. on $ \Gamma $
}
\end{array} \right\} \mbox{, for any $ \varepsilon \geq 0 $,}
\end{equation*}
and let us define \piegia{the projection} function $ \mathcal{T}_{B} : \R \to \overline{I_{B}} $, by putting:
\begin{equation}\label{T_B^rho}
	r \in \R \mapsto \mathcal{T}_{B} r := \left(r \vee ( \inf I_{B} ) \right) \wedge ( \sup I_{B} ) \in \overline{I_{B}}.
\end{equation}
Then, we easily check the following facts.
\begin{description}
\item[(Fact\,6)]If $ \varepsilon > 0 $, then the space $ \mathscr{V}_\varepsilon $ is a closed linear subspace 
in $ H^1(\Omega) \times H^1(\Gamma) $\piegia{. Otherwise  (i.e., if $ \varepsilon = 0 $),} 
the corresponding space $ \mathscr{V}_0 $ is a closed linear subspace in $ H^1(\Omega) \times H^{\frac{1}{2}}(\Gamma) $. 
Hence, when $ \varepsilon > 0 $ (resp.\ $ \varepsilon = 0 $), 
the space $ \mathscr{V}_\varepsilon $ (resp.\ $ \mathscr{V}_0 $) forms a Hilbert space 
endowed with the inner product in $ H^1(\Omega) \times H^1(\Gamma) $ 
(resp.\ $ H^1(\Omega) \times H^{\frac{1}{2}}(\Gamma) $). 
\item[(Fact\,7)] For any $ \varepsilon \geq 0 $, let us put: 
\begin{equation}\label{D-nu}
\D_{\varepsilon} := \left\{ \begin{array}{l|l}
[z, z_\Gamma] \in \mathscr{V}_\varepsilon & B(z) \in L^1(\Omega) \mbox{ and } B_\Gamma(z_\Gamma) \in L^1(\Gamma)
\end{array} \right\}.
\end{equation}
Then, the closures of $ \D_{\varepsilon} $, for $ \varepsilon \geq 0$, in the topology of $ \mathscr{H} $ 
coincide with the class $ \D_{*} $ given in~\eqref{D0},~i.e.
\begin{equation*}
\D_{*} = \overline{\D_{\varepsilon}} \mbox{ in $ \mathscr{H} $, for any $ \varepsilon \geq 0 $}.
\end{equation*}
\end{description}

Based on the above (A1)--(A3) and (Fact\,6)--(Fact\,7), the solutions to (ACE)$_{\varepsilon}$, for $ \varepsilon \geq 0 $, are defined as follows.
\begin{defn}[Definition of solutions]\label{Def.sol.}
\begin{em}
A pair $ [u, u_{\Gamma}] $ of functions $ u : [0,T] \to L^2(\Omega) $ and $ u_{\Gamma} : [0,T] \to L^2(\Gamma) $ is called a solution to (ACE)$_{\varepsilon}$, \piegia{iff} $ [u,u_{\Gamma}] $ fulfills the following conditions.
\begin{description}
\item[(S1)] $ [u,u_{\Gamma}] \in C([0,T];\mathscr{H}) \cap W_{\mathrm{loc}}^{1,2}((0,T];\mathscr{H}) \cap L^2(0,T;\mathscr{V}_\varepsilon) \cap L_{\mathrm{loc}}^{\infty}((0,T];\mathscr{V}_\varepsilon) $,
\\[1ex]
$ [u(0),u_{\Gamma}(0)] = [u_0, u_{\Gamma,0}] $ in $ \mathscr{H} $.

\item[(S2)] There exist functions $ \nu_{u} : Q \to \R^{N} $, $ \xi : (0,T) \to \piegia{L^2(\Omega)}$ and $ \xi_{\Gamma} : (0,T) \to \piegia{L^2(\Gamma)} $ such that:
\begin{equation*}
\begin{array}{l}
\nu_{u} \in L^{\infty}(Q) \mbox{ and } \nu_{u} \in \Sgn(\nabla u) \mbox{ a.e. in $ Q $, }
\\[2ex]
\xi \in L_{\mathrm{loc}}^2((0,T];L^2(\Omega)) \mbox{ and } \xi \in \beta(u) \mbox{ a.e. in $ Q $, }
\\[2ex]
\xi_{\Gamma} \in L_{\mathrm{loc}}^2((0,T];L^2(\Gamma)) \mbox{ and } \xi_{\Gamma} \in \beta_{\Gamma}(u_{\Gamma}) \mbox{ a.e. in $ \Sigma $, }
\end{array}
\end{equation*}
and
\begin{equation*}
\begin{array}{l}
\ds \int_{\Omega} \partial_{t}u(t) z \, dx + \int_{\Omega} (\nu_{u}(t) + \kappa^2 \nabla u(t)) \cdot \nabla z \, dx + \int_{\Omega} \bigl(\xi(t) + g(u(t)) \bigr) z \, dx
\\[2ex]
\ds +\int_{\Gamma} \partial_{t}u_{\Gamma}(t) z_{\Gamma} \, d\Gamma + \int_{\Gamma} \nabla_{\Gamma}(\varepsilon u_{\Gamma}(t)) \cdot \nabla_{\Gamma}(\varepsilon z_{\Gamma}) \, d\Gamma + \int_{\Gamma} \bigl( \xi_{\Gamma}(t) + g_{\Gamma}(u_{\Gamma}(t)) \bigr) z_{\Gamma} \, d\Gamma
\\[2ex]
\ds = \int_{\Omega} \theta(t)z \, dx + \int_{\Gamma} \theta_{\Gamma}(t) z_{\Gamma} \, d\Gamma \mbox{, for any $ [z,z_{\Gamma}] \in \mathscr{V}_\varepsilon $.}
\end{array}
\end{equation*}
\end{description}
\end{em}
\end{defn}

Now our Main Theorems are stated as follows.

\begin{mainTh}[well-posedness]\label{mainTh1}
\begin{em}
Let us assume (A1)--(A3) and let us fix \piegia{an} arbitrary $ \varepsilon \geq 0 $. Then, the following items hold.
\begin{description}
\item[(I-1)(Existence and uniqueness)] The system (ACE)$_{\varepsilon}$ admits a unique solution $ [u,u_{\Gamma}] $ and there exists a constant $ C_1 > 0 $, independent of the initial value $ [u_0,u_{\Gamma,0}] $ and the forcing term $ [\theta, \theta_{\Gamma}] $, such that:
\begin{equation}\label{est.1}
\begin{array}{lll}
\multicolumn{3}{l}{\ds\left| [u, u_{\Gamma}] \right|_{C([0,T];\mathscr{H})}^2 + \left| [\nabla u, \nabla_{\Gamma} (\varepsilon u_{\Gamma})] \right|_{L^2(0,T;\mathscr{H}^{N})}^2}
\\[2ex]
& & + \ds\left|\sqrt{t}\, [\partial_t u,\partial_t u_{\Gamma}] \right|_{L^2(0,T;\mathscr{H})}^2 + \sup_{t \in (0,T)} \left| \piegia{\sqrt{t}}\, [\nabla u(t), \nabla_{\Gamma} (\varepsilon u_{\Gamma}(t))] \right|_{\mathscr{H}^N}^2
\\[2ex]
\leq & \multicolumn{2}{l}{C_1 \left( 1 + \left| [u_0,u_{\Gamma,0}] \right|_{\mathscr{H}}^2 + \left| [\theta,\theta_{\Gamma}] \right|_{L^2(0,T;\mathscr{H})}^2 \right).}
\end{array}
\end{equation}
Moreover, if $ [u_0,u_{\Gamma,0}] \in \mathscr{D}_{\varepsilon} $, then there exists a constant $ C_2 > 0 $, independent of the initial value $ [u_0,u_{\Gamma,0}] $ and the forcing term $ [\theta,\theta_{\Gamma}] $, such that:
\begin{equation}\label{est.2}
\begin{array}{ll}
& \ds\left| [\partial_t u,\partial_t u_{\Gamma}] \right|_{L^2(0,T;\mathscr{H})}^2 + \sup_{t \in (0,T)} \left| [\nabla u(t), \nabla_{\Gamma} (\varepsilon u_{\Gamma}(t)) \right|_{\mathscr{H}^N}^2
\\[2ex]
\multicolumn{2}{l}{\leq C_2 \left( \begin{array}{rcl} 1 & + & \left| [u_0,u_{\Gamma,0}] \right|_{\mathscr{H}}^2 + \left| [\nabla u_0, \nabla_{\Gamma} (\varepsilon u_{\Gamma,0})] \right|_{\mathscr{H}^N}^2
\\[2ex]
& + & |B(u_0)|_{L^1(\Omega)} + |B_{\Gamma}(u_{\Gamma, 0})|_{L^1(\Gamma)}  + \left| [\theta,\theta_{\Gamma}] \right|_{L^2(0,T;\mathscr{H})}^2
\end{array} \right).}
\end{array}
\end{equation}

\item[(I-2)(Continuous-dependence)] For $ k = 1,2 $, let $ [u^{k},u_{\Gamma}^{k}] $ \piegia{denote} two solutions to \piegia{the problem}  (ACE)$_{\varepsilon}$ \piegia{corresponding to the} forcing pairs $ [\theta^{k},\theta_{\Gamma}^{k}] \in L^2(0,T;\mathscr{H}) $ and initial pairs $ [u_{0}^{k},u_{\Gamma,0}^{k}] \in \D_* $, respectively. Then, there exists a positive constant $ C_3 $, independent of the choices of $ [\theta^{k},\theta_{\Gamma}^{k}] $ and $ [u_{0}^{k},u_{\Gamma,0}^{k}] $, $ k=1,2 $, such that:
\begin{equation}\label{est.3}
\begin{array}{ll}
& \ds\left| [u^1 - u^2, u_{\Gamma}^1 - u_{\Gamma}^2] \right|_{C([0,T];\mathscr{H})}^2 
\piegia{{}+ \left| [\nabla (u^1 - u^2), \nabla_{\Gamma} (\varepsilon (u_{\Gamma}^1 - u_{\Gamma}^2))] \right|_{L^2(0,T;\mathscr{H}^{N})}^2}
\\[2ex]
\multicolumn{2}{l}{\ds\leq C_3 \left(\left| [u_0^1 - u_0^2, u_{\Gamma,0}^1 - u_{\Gamma,0}^2] \right|_{\mathscr{H}}^2 + \left| [\theta^1 - \theta^2, \theta_{\Gamma}^1 - \theta_{\Gamma}^2] \right|_{L^2(0,T;\mathscr{H})}^2 \right).}
\end{array}
\end{equation}
\end{description}
\end{em}
\end{mainTh}

\begin{mainTh}[$ \varepsilon $-dependence of solutions]\label{mainTh2}
\begin{em} 
Let $ \varepsilon_0 \geq 0 $ be a fixed constant. Let $ \{[\theta^{\varepsilon},\theta_{\Gamma}^{\varepsilon}]\}_{\varepsilon \ge 0} \subset L^2(0,T;\mathscr{H}) $ be a sequence of forcing pairs, let $ \{[u_0^{\varepsilon},u_{\Gamma,0}^{\varepsilon}] \in \D_{\varepsilon}\}_{\varepsilon \ge 0} \subset \mathscr{H} $ be a sequence of initial pairs and for any $ \varepsilon \geq 0 $, let $ [u^{\varepsilon},u_{\Gamma}^{\varepsilon}] $ be a solution to (ACE)$_{\varepsilon}$ corresponding to the forcing pair $ [\theta^{\varepsilon},\theta_{\Gamma}^{\varepsilon}] \in L^2(0,T;\mathscr{H}) $ and the initial pair $ [u_{0}^{\varepsilon},u_{\Gamma,0}^{\varepsilon}] \in \D_{\varepsilon} $. If:
\begin{equation}\label{MT2.0-0}
\left\{ \begin{array}{l} [\theta^{\varepsilon},\theta_{\Gamma}^{\varepsilon}] \to [\theta^{\varepsilon_0},\theta_{\Gamma}^{\varepsilon_0}] \mbox{ weakly in $ L^2(0,T;\mathscr{H}) $,}
\\[2ex]
[u_0^{\varepsilon},u_{\Gamma,0}^{\varepsilon}] \to [u_0^{\varepsilon_0},u_{\Gamma,0}^{\varepsilon_0}]  \mbox{ in $ \mathscr{H} $,}
\end{array} \right.
\mbox{ as $ \varepsilon \to \varepsilon_0 $,}
\end{equation}
then:
\begin{equation}\label{MT2.0-1}
\begin{array}{ll}
[u^{\varepsilon},u_{\Gamma}^{\varepsilon}] \to [u^{\varepsilon_0},u_{\Gamma}^{\varepsilon_0}] & \mbox{ in $ C([0,T];\mathscr{H}) $ 
and}
\\[0.5ex]
                                                                                                                                 & \mbox{ in $ L^2(0,T;\mathscr{V}_{0}) $, as $ \varepsilon \to \varepsilon_0 $.}
\end{array}
\end{equation}
In particular, if $ \varepsilon_0 > 0 $, then:
\begin{equation}\label{MT2.0-2}
u_{\Gamma}^\varepsilon \to u_{\Gamma}^{\varepsilon_0} \mbox{ in $ L^2(0, T; H^1(\Gamma)) $ as $ \varepsilon \to \varepsilon_0 $.}
\end{equation}
\end{em}
\end{mainTh}

\section{Key-Lemmas}
\ \ \vspace{-3ex}

In this Section, we specify the essential points in the proofs of  Main Theorems, in forms of Key-Lemmas.

In any case of $ \varepsilon \geq 0 $, the keypoint will be to reformulate the system (ACE)$_\varepsilon$ \piegia{as} the following Cauchy problem \piegia{(CP)$_\varepsilon$ for an evolution equation:}
\medskip

\noindent
(CP)$_\varepsilon$:
\begin{equation*}
\left\{ ~ \parbox{10cm}{
$ U'(t) +\partial \Phi_\varepsilon(U(t)) +\mathcal{G}(U(t)) \ni \Theta(t) $ in $ \mathscr{H} $, a.e. $ t \in (0, T) $,
\\[1ex]
$ U(0) = U_0 $ in $ \mathscr{H} $.
} \right.
\end{equation*}
In the context:
\begin{description}
\item[--]the unknown $ U \in C([0, T]; \mathscr{H}) $ corresponds to the solution pair $ [u, u_\Gamma] $ of the system (ACE)$_\varepsilon$,~i.e.,
$ U(t) = [u(t), u_\Gamma(t)] $ in $ \mathscr{H} $, for any $ t \in [0, T] $ with the initial pair $ U_0 = [u_0, u_{\Gamma, 0}] $ in $ \mathscr{H} $;
\item[--]$ \partial \Phi_\varepsilon $ denotes the subdifferential of a proper l.s.c.\ and convex function 
$ \Phi_\varepsilon : \mathscr{H} \to [0, \infty] $, defined as
\end{description}
\begin{equation}\label{Phi_nu}
\begin{array}{rcl}
U & = & [u, u_\Gamma] \in \mathscr{H} \mapsto \Phi_\varepsilon(U) = \Phi_\varepsilon(u, u_\Gamma)
\\[2ex]
& := & \left\{ \begin{array}{lll}
\multicolumn{3}{l}{\ds \int_\Omega \left( |\nabla u| +\frac{\kappa^2}{2} |\nabla u|^2 \right) \, dx +\int_\Omega B(u) \, dx}
\\[2ex]
&& \ds + \frac{1}{2} \int_\Gamma |\nabla_\Gamma(\varepsilon u_\Gamma)|^2 \, d \Gamma +\int_\Gamma B_\Gamma(u_\Gamma) \, d \Gamma,
\\[2ex]
& \multicolumn{2}{l}{\mbox{if $ U = [u, u_\Gamma] \in \mathscr{D}_\varepsilon $,}}
\\[2ex]
\infty, & \multicolumn{2}{l}{\mbox{otherwise,}}
\end{array} \right.
\ \\
\\[-2ex]
\multicolumn{3}{l}{\mbox{for $ \varepsilon \geq 0 $;}}
\end{array}
\end{equation}
\begin{description}
\item[--]$ \mathcal{G} : \mathscr{H} \to \mathscr{H} $ is a \piegia{Lipschitz continuous} operator, defined as
\begin{equation*}
U = [u, u_\Gamma] \in \mathscr{H} \mapsto \mathcal{G}(U) = \mathcal{G}(u, u_\Gamma) := [g(u), g_\Gamma(u_\Gamma)] \in \mathscr{H}\piegia{,}
\end{equation*}
\piegia{where $g$ and $g_\Gamma$ can be meant as Lipschitz continuous extensions outside $\overline{I_B}$ of the functions $g$ and $g_\Gamma$ defined in (A2);}
\item[--]the forcing term $ \Theta $ corresponds to the forcing pair $ [\theta, \theta_\Gamma] $ of (ACE)$_\varepsilon$, i.e., 
$ \Theta = [\theta, \theta_\Gamma] $ in  $ L^2(0, T; \mathscr{H}) $.
\end{description}

\begin{rem}\label{Rem.ConvEs}
\begin{em}
	For any $ \varepsilon \geq 0 $, we can see that the convex function $ \Phi_\varepsilon $, \piegia{defined in} \eqref{Phi_nu}, corresponds to the convex part of the free-energy, given in \eqref{free-energy}. In addition, the subdifferentials $ \partial \Phi_{\varepsilon} $, for $ \varepsilon \geq 0 $, are maximal monotone graphs in $ \mathscr{H} \times \mathscr{H} $. So, the well-posedness for the Cauchy problem (CP)$_{\varepsilon}$ will be verified, immediately, by applying general theories \piegia{for evolution equations, e.g.,} \cite{Barbu, Brezis, I-Y-K}. 
\end{em}
\end{rem}

In the light of Remark \ref{Rem.ConvEs}, the essential points in Main Theorem \ref{mainTh1} will be to show a certain association between our system (ACE)$_{\varepsilon}$ and the Cauchy problem (CP)$_{\varepsilon}$ for any $ \varepsilon \geq 0 $. To this end, we need to prepare a class of relaxed convex functions 
$$ \piegia{\{ \Phi_{\varepsilon, \delta}^{\lambda} \ |\ \varepsilon \geq 0,\  0 < \delta,\, \lambda \leq 1 \} ,} $$ defined as follows:
\begin{equation}\label{ap.Phi}%\label{Phi_nu,dlt^lmd}
\begin{array}{rcl}
U & = & [u, u_\Gamma] \in \mathscr{H} \mapsto \Phi_{\varepsilon, \delta}^\lambda(U) = \Phi_{\varepsilon, \delta}^\lambda(u, u_\Gamma)
\\[2ex]
& := & \left\{ \begin{array}{lll}
\multicolumn{3}{l}{\ds \int_\Omega \left( f_\delta(\nabla u) +\frac{\kappa^2}{2} |\nabla u|^2 \right) \, dx +\int_\Omega B^\lambda(u) \, dx}
\\[2ex]
&& \ds +\frac{1}{2} \int_\Gamma |\nabla_\Gamma(\varepsilon u_\Gamma)|^2 \, d \Gamma +\int_\Gamma B_\Gamma^\lambda(u_\Gamma) \, d \Gamma,
\\[2ex]
& \multicolumn{2}{l}{\mbox{if $ U = [u, u_\Gamma] \in \mathscr{V}_\varepsilon $,}}
\\[2ex]
\infty, & \multicolumn{2}{l}{\mbox{otherwise,}}
\end{array} \right.
\ \\
\\[-2ex]
\multicolumn{3}{l}{\mbox{for $ \varepsilon \geq 0 $ and $ 0 < \delta, \/ \lambda \leq 1 $.}}
\end{array}
\end{equation}
In the context, $ \{ f_\delta \}_{0 < \delta \leq 1} $, $ \{ B^\lambda \}_{0 < \lambda \leq 1} $ and $ \{ B_\Gamma^\lambda \}_{0 < \lambda \leq 1} $ are sequences of functions, prescribed under the following assumptions.
\begin{description}
\item[(A4)]$ \{ f_\delta \}_{ 0 < \delta \leq 1} \subset C^1(\R^N) $ is a sequence of convex functions and $C^1$-regularizations for the Euclidean norm $ |{}\cdot{}| \in W^{1, \infty}(\R^N) $, such that:
\begin{equation*}
f_{\delta}(0) = 0 \mbox{ and } f_{\delta}(\omega) \geq 0 \mbox{, for any $ \omega \in \R^N $ and any $ 0 < \delta \leq 1 $,}
\end{equation*}
\begin{equation*}
\left\{\begin{array}{l}
f_{\delta}(\omega) \to |{}\omega{}| \mbox{, for any $ \omega \in \R^N $,}
\\[2ex]
f_{\delta} \to |{}\cdot{}| \mbox{ on $ \R^N $, in the sense of \piegia{Mosco},}
\end{array}\right.
\mbox{as $ \delta \downarrow 0 $,}
\end{equation*}
and there exists a $\delta$-independent constant $ C_0 > 0 $, satisfying:
\begin{equation*}
|{}\nabla f_{\delta} (\omega){}| \leq C_0(|\omega| + 1) \mbox{, for any $ 0 < \delta \leq 1 $ and $ \omega \in \R^N $.}
\end{equation*}

\item[(A5)]$ \{ B^\lambda \}_{0 < \lambda \leq 1} $ and $ \{ B_\Gamma^\lambda \}_{0 < \lambda \leq 1} $ 
are sequences of Moreau--Yosida regularizations of \piegia{the} convex functions $ B $ and $ B_\Gamma $, respectively, i.e., 
$ \{ B^\lambda \}_{0 < \lambda \leq 1} \subset C^1(\R) $, $ \{ B_\Gamma^\lambda \}_{0 < \lambda \leq 1} \subset C^1(\R) $, and
\begin{equation*}
\left\{ 
\begin{array}{l}
\tau \in \R \mapsto B^\lambda(\tau) := \inf \left\{ \begin{array}{l|l}
\ds \frac{1}{2 \lambda} |\tilde{\tau} -\tau|^2 + B(\tilde{\tau}) & \tilde{\tau} \in \R
\end{array} \right\}, 
\\[2ex]
\tau \in \R \mapsto B_\Gamma^\lambda(\tau) := \inf \left\{ \begin{array}{l|l}
\ds \frac{1}{2 \lambda} |\tilde{\tau} - \tau|^2 + B_\Gamma(\tilde{\tau}) & \tilde{\tau} \in \R
\end{array} \right\}, 
\end{array} \right. \mbox{ for any $ 0 < \lambda \leq 1 $.}
\end{equation*}
\end{description}

Now, the first Key-Lemma is concerned with the representations of the subdifferentials $ \partial \Phi_{\varepsilon,\delta}^{\lambda} \subset \mathscr{H} \times \mathscr{H} $ of the relaxed convex functions $ \Phi_{\varepsilon,\delta}^{\lambda} $, for $ \varepsilon \geq 0 $, $ 0 < \delta, \/ \lambda \leq 1 $.
\begin{keyLem}\label{keyLem0}
Let us fix $ \varepsilon \ge 0 $, $ 0 < \delta, \lambda \le 1 $. Let us put:
\begin{equation*}
\mathscr{D}_{\varepsilon, \delta}^{\lambda} :=\left\{ \begin{array}{l|l} [u,u_{\Gamma}] \in \mathscr{V}_{\varepsilon}  & \begin{array}{l} \nabla f_{\delta} (\nabla u) + \kappa^2 \nabla u \in \bm{L}_{\rm div}^2(\Omega) \\[1ex] -{\mit \Delta}_{\Gamma} (\varepsilon^2 u_{\Gamma}) + \left[(\nabla f_{\delta} (\nabla u) + \kappa^2 \nabla u) \cdot n_{\Gamma} \right]_{_{\Gamma}} \in L^2(\Gamma) \end{array} \end{array} \right\},
\end{equation*}
and let us define an operator $ \mathcal{A}_{\varepsilon, \delta}^{\lambda} : \mathscr{D}_{\varepsilon,\delta}^{\lambda} \subset \mathscr{H} \to \mathscr{H} $, by putting:
\begin{equation*}
\begin{array}{rcl}
[u, u_{\Gamma}] \in \mathscr{D}_{\varepsilon,\delta}^{\lambda} & \mapsto & \mathcal{A}_{\varepsilon,\delta}^{\lambda}[u,u_{\Gamma}] \\[1ex]
							       & := & {\rule{0pt}{18pt}}^{t} \hspace{-0.5ex} \left[
\begin{array}{l} 
-\mathrm{div}(\nabla f_{\delta} (\nabla u) + \kappa^2 \nabla u) + \beta^{\lambda}(u) \\[1ex]
-{\mit \Delta}_{\Gamma} (\varepsilon^2 u_{\Gamma}) + [(\nabla f_{\delta} (\nabla u) + \kappa^2 \nabla u) \cdot n_{\Gamma} ]_{_{\Gamma}} + \beta_{\Gamma}^{\lambda}(u_{\Gamma})
\end{array}\right] \mbox{ in $ \mathscr{H} $.}
\end{array}
\end{equation*}
Then, $ \partial \Phi_{\varepsilon,\delta}^{\lambda} = \mathcal{A}_{\varepsilon,\delta}^{\lambda} $ in $ \mathscr{H} \times \mathscr{H} $,~i.e.
\begin{equation*}
\begin{array}{c}
D(\partial \Phi_{\varepsilon,\delta}^{\lambda}) = \mathscr{D}_{\varepsilon,\delta}^{\lambda} \mbox{ and } \partial \Phi_{\varepsilon,\delta}^{\lambda} (u,u_{\Gamma}) = \mathcal{A}_{\varepsilon,\delta}^{\lambda}[u,u_{\Gamma}] \mbox{ in $ \mathscr{H} $, }\\[1ex]
\mbox{ for any $ [u,u_{\Gamma}] \in D(\partial \Phi_{\varepsilon,\delta}^{\lambda}) $.}
\end{array}
\end{equation*}
\end{keyLem} 

The second Key-Lemma is concerned with the continuous dependence between the convex functions $ \Phi_{\varepsilon} $ for $ \varepsilon \geq 0 $, and the relaxations of those.
\begin{keyLem}[Continuous dependence of the convex functions]\label{keyLem2}
Let 
$$ \piegia{\{ \varepsilon_n \}_{n = 0}^\infty \subset [0, \infty) ,\quad  \{ \delta_n \}_{n = 1}^\infty \subset (0, 1]  \quad \hbox{and} \quad \{ \lambda_n \}_{n = 1}^\infty \subset (0, 1] }$$ 
be arbitrary sequences such that:
\begin{equation*}
\varepsilon_n \to \varepsilon_0, \quad \delta_n \downarrow 0 \ \mbox{ and } \ \lambda_n \downarrow 0 \mbox{, \  as $ n \to \infty $.}
\end{equation*}
Then, for the sequence $ \{ \Phi_{\varepsilon_n, \delta_n}^{\lambda_n} \}_{n = 1}^\infty $ of convex functions, it holds that:
\begin{equation*}
\Phi_{\varepsilon_n, \delta_n}^{\lambda_n} \to \Phi_{\varepsilon_0} \mbox{ on $ \mathscr{H} $, in the sense of  \piegia{Mosco}, as $ n \to \infty $.}
\end{equation*}
\end{keyLem}

On the basis of Key-Lemmas \ref{keyLem0}--\ref{keyLem2}, we prove the third Key-Lemma, concerned with representations of the subdifferentials $ \partial \Phi_{\varepsilon} \subset \mathscr{H} \times \mathscr{H} $ of $ \Phi_{\varepsilon} $, for $ \varepsilon \geq 0 $.
\begin{keyLem}\label{keyLem1}
For any $ \varepsilon \geq 0 $, the following two items are equivalent.
\begin{description}
\item[(Key\,0)]$ U = [u, u_\Gamma] \in D(\partial \Phi_\varepsilon) $ and $ U^* = [u^*, u_\Gamma^*] \in \partial \Phi_\varepsilon(U) = \partial \Phi_\varepsilon(u, u_\Gamma) $ in $ \mathscr{H} $.
\item[(Key\,1)]$ U = [u, u_\Gamma] \in \mathscr{D}_\varepsilon $, and there exists $ \nu_u \in L^\infty(\Omega)^N $ and $ [\xi, \xi_\Gamma] \in \mathscr{H} $, such that:
\begin{equation}\label{key1.1}
\left\{ ~ \parbox{9cm}{
$ \nu_u \in {\rm Sgn}(\nabla u) $ and $ \xi \in \beta(u) $, a.e. in $ \Omega $,
\\[1ex]
$ \xi_\Gamma \in \beta_\Gamma(u_\Gamma) $, a.e. on $ \Gamma $,
} \right.
\end{equation}
\begin{equation}\label{key1.2}
\left\{ ~ \parbox{9cm}{
$ \nu_u +\kappa^2 \nabla u \in \bm{L}_{\rm div}^2(\Omega) $,
\\[1ex]
$ -{\mit \Delta}_\Gamma (\varepsilon^2 u_\Gamma) +[(\nu_u +\kappa^2 \nabla u) \cdot n_\Gamma]_{_{\Gamma}} \in L^2(\Gamma) $,
} \right.
\end{equation}
and
\begin{equation}\label{key1.3}
\left\{ ~ \parbox{9cm}{
$ u^*  = -{\rm div} \, \bigl( \nu_u  +\kappa^2 \nabla u \bigr) + \xi $ in $ L^2(\Omega) $,
\\[1ex]
$ u_\Gamma^* = -{\mit \Delta}_\Gamma (\varepsilon^2 u_\Gamma) +[(\nu_u +\kappa^2 \nabla u) \cdot n_\Gamma]_{_{\Gamma}} + \xi_{\Gamma} $ in $ L^2(\Gamma) $.
} \right.
\end{equation}
\end{description}
\end{keyLem}

The last Key-Lemma \ref{keyLem1} \piegia{is useful} to guarantee the association between (ACE)$_{\varepsilon}$ and (CP)$_{\varepsilon}$ for $ \varepsilon \geq 0 $, via the representations of subdifferentials.

%%%%%

\section{Proofs of Key-Lemmas}
\ \ \vspace{-3ex}

In this section, we prove three Key-Lemmas for our Main Theorems. To this end, we first prepare the following lemma.
\begin{lem}\label{auxLem}
Let $ (S, \mathcal{B}, \mu) $ be a measure space with a $ \sigma $-algebra $ \mathcal{B} $ and a finite Radon measure $ \mu $. 
Let $ X $ be a (real) Hilbert space.  
\piegia{Let $ \Psi : X \rightarrow (-\infty, \infty] $ be a proper l.s.c.\ and convex function, 
and let $ \{ \Psi_n \}_{n = 1}^\infty $ be a sequence of proper l.s.c.\ and convex functions 
$ \Psi_n : X \rightarrow (-\infty, \infty] $, $ n \in \N $,  such that:}
\begin{equation}\label{axLem01}
\Psi_n \to \Psi \mbox{ on $ X $, in the sense of \piegia{Mosco}, as $ n \to \infty $.}
\end{equation}
Then, the following two items hold.
\begin{description}
\item[(I)] There exist two constants $ c_0,d_0 > 0 $, independent of $ n $, such that:
\begin{equation}\label{axLem02}
\Psi_n(z) + c_0|z|_X + d_0 \geq 0 \mbox{, for any $ z \in X $ and any $ n \in \N $.}
\end{equation}
\item[(I\hspace{-.1ex}I)] \piegia{The} sequence $ \{ \hat{\Psi}_n \}_{n = 1}^{\infty} $ of proper l.s.c.\ and convex functions
\begin{equation}\label{axLem03}
\begin{array}{rcl}
\zeta \in L^2(S;X) & \mapsto & \hat{\Psi}_n(\zeta) := \left\{\begin{array}{ll} {\ds \int_{S} \Psi_n(\zeta) \, d\mu \mbox{, if $ \Psi_n(\zeta) \in L^1(S) $, }}
\\[0.5ex]
\infty \mbox{, otherwise, }
\end{array}\right.
\\[4ex]
&&  \mbox{for $ n = 1, 2, 3, \dots, $}
\end{array}
\end{equation}
converges to \piegia{the convex function}
\begin{equation}\label{axLem04}
\begin{array}{rcl}
\zeta \in L^2(S;X) & \mapsto & \hat{\Psi}(\zeta) := \left\{\begin{array}{ll} {\ds \int_{S} \Psi(\zeta) \, d\mu \mbox{, if $ \Psi(\zeta) \in L^1(S) $, }}
\\[0.5ex]
\infty \mbox{, otherwise,}
\end{array}\right.
\end{array}
\end{equation}
on the Hilbert space $ L^2(S;X) $, in the sense of  \piegia{Mosco}, as $ n \to \infty $.
\end{description}
\end{lem}
\textbf{Proof.} This lemma can be proved by means of similar demonstration techniques as in \cite[Appendix]{GKY}.
\piegia{However, we report the proof for the reader's convenience.}

First, we show the item (I). To this end, let us assume that:
\begin{equation}\label{axLem05}
\begin{array}{c}
\Psi_{n_k}(y_k) + k^2(|y_k|_{X} + 1) < 0,\\[1ex]
\mbox{ for some $ \{ n_k \}_{ k = 1 }^{\infty} \subset \{ n \}_{n = 1}^{\infty} $, $ \{ y_k \}_{k = 1}^{\infty} \subset X $}
\end{array}
\end{equation}
to derive a contradiction.

Let us fix any $ z_0 \in D(\Psi) $. Then, by \eqref{axLem01}, we find a sequence $ \{ \hat{z}_n \}_{n = 1}^{\infty} \subset X $, such that:
\begin{equation}\label{axLem06}
\hat{z}_n \to z_0 \mbox{ in $ X $ and } \Psi_n(\hat{z}_n) \to \Psi(z_0) \mbox{ as $ n \to \infty $. }
\end{equation}
Here, we define:
\begin{equation}\label{axLem07}
z_k := \varepsilon_k y_k + ( 1 - \varepsilon_k ) \hat{z}_{n_k} \mbox{ in $ X $, for $ k =1,2,3,\dots $,}
\end{equation}
with
\begin{equation}\label{axLem08}
\varepsilon_k := \frac{1}{k(1 + |y_k|_X)} \in (0,1) \mbox{, for $ k=1,2,3,\dots $. }
\end{equation}
Then, it follows from \eqref{axLem05}--\eqref{axLem08} that \piegia{$|\varepsilon_k y_k|\leq 1/k$, whence:}
\begin{equation*}
z_{k} \to z_0 \mbox{ in $ X $ as $ k \to \infty $, }
\end{equation*}
and subsequently, it follows from \eqref{axLem01} that:
\begin{equation}\label{axLem09}
\varlimsup_{k \to \infty} \Psi_{n_k}(z_k) \geq \varliminf_{k \to \infty} \Psi_{n_k}(z_k) \geq \Psi(z_0).
\end{equation}

In the meantime, in the light of \eqref{axLem07}--\eqref{axLem08}, and the convexity of $ \Psi_{n_k} $, for $ k \in \N $,
\begin{align*}
\varlimsup_{ k \to \infty } \Psi_{n_k}(z_k) & \leq \varlimsup_{k \to \infty} \varepsilon_k \Psi_{n_k}(y_k) + \lim_{k \to \infty} (1 - \varepsilon_k) \Psi_{n_k} (\hat{z}_{n_k})
\\[1ex]
& \leq \varlimsup_{k \to \infty} \left( \varepsilon_k \big({- k^2}(1 + |y_k|_X)\big)\right) \piegia{{}+\Psi(z_0)} 
\\[1ex]
& =  - \lim_{k \to \infty} k  \piegia{{}+\Psi(z_0) = - \infty.}
\end{align*}
This contradicts with \eqref{axLem09}.

Next, we show the item (I\hspace{-.1ex}I). According to \cite[Theorem 3.26]{Attouch}, it is sufficient (equivalent) to check the following two conditions:
\begin{description}
\item[(ii-1)] $ \zeta_n := (\mathcal{I}_{L^2(S;X)} + \lambda \partial \hat{\Psi}_n)^{-1} \xi \to \zeta := (\mathcal{I}_{L^2(S;X)} + \lambda \partial \hat{\Psi})^{-1} \xi \mbox{ in $ L^2(S;X) $ as $ n \to \infty $,} $
\\[1ex]
\mbox{ for any $ \lambda > 0 $ and any $ \xi \in L^2(S;X) $};
\item[(ii-2)] there exists $ [\zeta, \eta] \in \partial \hat{\Psi} $ in $ L^2(S;X) \times L^2(S;X) $ and a sequence $ \{ [\zeta_n, \eta_n] \in \hat{\Psi}_n\}_{n=1}^{\infty} \subset L^2(S;X) \times L^2(S;X) $ such that $ [\zeta_n,\eta_n] \to [\zeta,\eta] $ in $ L^2(S;X) \times L^2(S;X) $ and $ \hat{\Psi}_n(\zeta_n) \to \hat{\Psi}(\zeta) $ as $ n \to \infty $.
\end{description}

For the verification of (ii-1), let us fix any $ \lambda > 0 $ and any $ \xi \in L^2(S;X) $. Then\piegia{,} invoking 
\piegia{\cite[Proposition~2.16]{Brezis}, \cite[Theorem~3.26]{Attouch},} \eqref{axLem01} and 
\eqref{axLem03}--\eqref{axLem04}, we infer that:
\begin{equation}\label{axLem12}
\left\{\begin{array}{l}
(\xi - \zeta)(\sigma) \in \lambda \partial \piegia{\Psi}(\zeta(\sigma)) \mbox{ in $ X $, }
\\[0.5ex]
(\xi - \zeta_n)(\sigma) \in \lambda \partial \piegia{\Psi}_n(\zeta_n(\sigma)) \mbox{ in $ X $, $ n = 1,2,3,\dots $,}
\end{array}\right.
\mbox{ for $ \mu $-a.e. $ \sigma \in S $,}
\end{equation}
and
\begin{equation}\label{axLem13}
\zeta_n(\sigma) \to \zeta(\sigma) \mbox{ in $ X $ as $ n \to \infty $, for $ \mu $-a.e. $ \sigma \in S $.}
\end{equation}
Also, by using the sequence $ \{ \hat{z}_n \}_{n = 1}^{\infty} \subset X $ as in \eqref{axLem06}, it is seen that:
\begin{equation}\label{axLem14}
\begin{array}{c}
\bigl( ( \xi - \zeta_n)(\sigma), \zeta_n(\sigma) - \hat{z}_n \bigr)_X \geq \lambda \hat{\Psi}_n(\zeta_n(\sigma)) - \lambda \hat{\Psi}_n(\hat{z}_n),
\\[1ex]
\mbox{ for any $ n \in \N $ and $ \mu $-a.e. $ \sigma \in S $}.
\end{array}
\end{equation}
Additionally, by virtue of the item (I), \eqref{axLem06}, \eqref{axLem14} and \piegia{the Schwarz and Young inequalities}, we can compute that:
\begin{equation*}
\begin{array}{rcl}
|\zeta_n(\sigma)|_X^2 & \leq & \ds(\zeta_n(\sigma),\hat{z}_n)_X + (\xi(\sigma),\zeta_n(\sigma) - \hat{z}_n)_X + \lambda \hat{\Psi}_n(\hat{z}_n) - \lambda \hat{\Psi}_n(\zeta_n(\sigma))
\\[2ex]
& \leq & \ds|\zeta_n(\sigma)|_X|\hat{z}_n|_X + |\zeta_n(\sigma)|_X|\xi(\sigma)|_X + |\xi(\sigma)|_X|\hat{z}_n|_X
\\[2ex]
& & \ds+ \lambda \hat{\Psi}_n(\hat{z}_n) + \lambda \bigl( c_0|\zeta_n(\sigma)|_X + d_0 \bigr)
\\[2ex]
& \leq & \ds\frac{3}{4} |\zeta_n(\sigma)|_X^2 + \frac{5}{4} |\xi(\sigma)|_X^2 + \lambda^2 c_0^2 + \lambda d_0 + \lambda \hat{\Psi}_n(\hat{z}_n) + 2|\hat{z}_n|_X^2,
\end{array}
\end{equation*}
and therefore,
\begin{equation}\label{axLem15}
|\zeta_n(\sigma)|_X^2 \leq 5|\xi(\sigma)|_X^2 + \hat{M}_1 \mbox{, for $ \mu $-a.e. $ \sigma \in S $ and any $ n \in \N $,}
\end{equation}
where
\begin{equation*}
\hat{M}_1:= 4\bigl( \lambda (c_0 + d_0) + 1 \bigr)^2 + 4 \sup_{n \in \N} (\lambda \hat{\Psi}_n(\hat{z}_n) + 2|\hat{z}_n|_X^2).
\end{equation*}
In view of these, the condition (ii-1) will be obtained as a consequence of \eqref{axLem13}, \eqref{axLem15} and Lebesgue's dominated convergence theorem.

Finally, for the verification of (ii-2), \piegia{we} consider the class of functions 
$ \{\zeta,\zeta_n | n \in \N \} \subset L^2(S;X) $ as in (ii-1) with fixed $ \lambda > 0 $ and $ \xi \in L^2(S;X) $, and let us set:
\begin{equation}\label{axLem16}
\left\{\begin{array}{l}
\ds \eta := \ds\frac{\xi - \zeta}{\lambda} \mbox{ in $ L^2(S;X) $,}
\\[1.5ex]
\ds \eta_n := \ds\frac{\xi - \zeta_n}{\lambda} \mbox{ in $ L^2(S;X) $, $ n=1,2,3, \dots $.}
\end{array} \right.
\end{equation}
Also, let us denote by $ \Psi^{\lambda} : X \to \R $, $ \hat{\Psi}^{\lambda} : L^2(S;X) \to \R $, $ \Psi_n^{\lambda} : X \to \R $ and $ \hat{\Psi}_n^{\lambda} : L^2(S;X) \to \R $, $ n \in \N $, the Moreau-Yosida regularizations of convex functions $ \Psi $, $ \hat{\Psi} $, $ \Psi_n $ and $ \hat{\Psi}_n $, $ n \in \N $, respectively. Then, by \cite[\piegia{Theorem~2.9, p.~48}]{Barbu}, \cite[Proposition~2.11]{Brezis}, (ii-1) and \eqref{axLem16}, we immediately have:
\begin{equation}\label{axLem17}
\left\{\begin{array}{l}
\ds \eta = \ds\partial \hat{\Psi}^{\lambda}(\xi) \in \partial \hat{\Psi}(\zeta) \mbox{ in $ L^2(S;X) $,}
\\[0.5ex]
\ds \eta_n = \ds\partial \hat{\Psi}_n^{\lambda}(\xi) \in \partial \hat{\Psi}_n(\zeta_n) \mbox{ in $ L^2(S;X) $, $ n=1,2,3,\dots $,}
\end{array} \right.
\end{equation}
\begin{equation}\label{axLem18}
\eta_n = \ds{\partial \hat{\Psi}_n^{\lambda}(\xi) = \frac{\xi - \zeta_n}{\lambda}} \to \ds{\eta = \partial \hat{\Psi}^{\lambda}(\xi) = \frac{\xi - \zeta}{\lambda} \mbox{ in $ L^2(S;X) $ as $ n \to \infty $.}}
\end{equation}
In particular, \cite[Proposition 2.16]{Brezis} and \eqref{axLem18} \piegia{enable} us to say that:
\begin{equation}\label{axLem19}
\partial \Psi_n^{\lambda} (\xi(\sigma)) \to \partial \Psi^{\lambda}(\xi(\sigma)) \mbox{ in $ X $ as $ n \to \infty $, for $ \mu $-a.e. $ \sigma \in S $,}
\end{equation}
by taking a subsequence if necessary. Besides, from \cite[Theorem 3.26]{Attouch} and \eqref{axLem01}, it follows that:
\begin{equation}\label{axLem20}
\Psi_n^{\lambda}(\xi(\sigma)) \to \Psi^{\lambda} (\xi(\sigma)) \mbox{ as $ n \to \infty $, for $ \mu $-a.e. $ \sigma \in S $.}
\end{equation}
On account of \cite[\piegia{Theorem~2.9, p.~48}]{Barbu}, \cite[Proposition 2.11]{Brezis} and \eqref{axLem19}--\eqref{axLem20}, it is inferred that:
\begin{equation}\label{axLem21}
\begin{array}{rcl}
\ds\Psi_n(\zeta_n(\sigma)) & = & \ds\Psi_n^{\lambda}(\xi(\sigma)) - \frac{1}{2 \lambda} | (\xi - \zeta_n)(\sigma)|_X^2 \\[2ex]
                                      & \to & \ds \Psi^{\lambda}(\xi(\sigma)) - \frac{1}{2 \lambda} | (\xi - \zeta)(\sigma)|_X^2 \\[2ex]
                                      & = & \ds \Psi(\zeta(\sigma)) \mbox{ as $ n \to \infty $, for $ \mu $-a.e. $ \sigma \in S $.}                                      
\end{array}
\end{equation}
\ryoken{Furthermore, invoking \eqref{axLem14}--\eqref{axLem15}, the item(I), and using 
the sequence $ \{ \hat{z}_n \}_{n=1}^{\infty} $ as in \eqref{axLem06} and \giapie{the 
Schwarz and Young inequalities}, we obtain that:}
\ryoken{
\begin{equation}\label{axLem22}
\begin{array}{rcl}
|\Psi_n(\zeta_n(\sigma))| & \leq & \Psi_n(\zeta_n(\sigma)) \vee (c_0 | \zeta_n(\sigma) |_X + d_0)
\\[2ex]
& \leq & \ds| \Psi_n(\hat{z}_n) | + \left| \left(\frac{(\xi - \zeta_n)(\sigma)}{\lambda}, \hat{z}_n - \zeta_n(\sigma)\right)_X \right| + c_0|\zeta_n(\sigma)|_X + d_0
\\[2ex]
& \leq & \ds| \Psi_n(\hat{z}_n) | + \frac{1}{\lambda}|\xi(\sigma)|_X|\hat{z}_n|_X + \frac{1}{\lambda}|\xi(\sigma)|_X|\zeta_n(\sigma)|_X 
\\[2ex]
& & \qquad \ds+ \frac{1}{\lambda}|\zeta_n(\sigma)|_X|\hat{z}_n|_X + \frac{1}{\lambda}|\zeta_n(\sigma)|_X^2 + c_0|\zeta_n(\sigma)|_X + d_0
\\[2ex]
& \leq & \ds|\Psi_n(\hat{z}_n)| + \frac{3}{\lambda}|\zeta_n(\sigma)|_X^2 + \frac{1}{\lambda}|\xi(\sigma)|_X^2 + \frac{1}{\lambda}|\hat{z}_n|_X^2 + \frac{\lambda}{4}c_0^2 + d_0
\\[2ex]
& \leq & \ds\frac{3}{\lambda}\left( 5|\xi(\sigma)|_X^2 + \hat{M}_1 \right) + \frac{1}{\lambda}|\xi(\sigma)|_X^2 
\\[2ex]
& & \qquad \ds+ |\Psi_n(\hat{z}_n)| + \frac{1}{\lambda}|\hat{z}_n|_X^2 + \frac{\lambda}{4}c_0^2 + d_0
\\[2ex]
& \leq & \ds\frac{16}{\lambda}|\xi(\sigma)|_X^2 + \hat{M}_2 \mbox{, for $ \mu $-a.e. $ \sigma \in S $ and any $ n \in \N $,}
\end{array}
\end{equation}
where
\begin{equation*}
\hat{M}_2 := \frac{3}{\lambda} \hat{M}_1 + \frac{\lambda}{4}c_0^2 + d_0 + \sup_{n \in \N} \left( \Psi_n(\hat{z}_n) + \frac{1}{\lambda}|\hat{z}_n|_X^2 \right).
\end{equation*}
}%
%%%%%%%%%%%%%%
\begin{comment}
Furthermore, invoking \eqref{axLem14} \piegia{{\bf (INSTEAD OF \eqref{axLem15} AS IT WAS IN YOUR VERSION. MOREOVER, PLEASE ADD SOME DETAILS TO THE INEQUALITIES BELOW IN ORDER THE READERS CAN UNDERSTAND THE VALUE OF CONSTANTS)}} and the item (I), and using the sequence $ \{ \hat{z}_n \}_{n=1}^{\infty} $ as in \eqref{axLem06}, we observe that:
\begin{equation}\label{axLem22}
\begin{array}{rcl}
| \Psi_n(\zeta_n(\sigma)) | & \leq & \ds\Psi_n(\zeta_n(\sigma)) \vee (c_0 |\zeta_n(\sigma)|_X + d_0)
\\[2ex]
                      & \leq & \ds |\Psi_n(\hat{z}_n)| + \left| \left( \frac{(\xi - \zeta_n)(\sigma)}{\lambda}, \hat{z}_n - \zeta_n(\sigma) \right)_X \right| + c_0 |\zeta_n(\sigma)|_X + d_0
\\[2ex]
                      & \leq & \ds \piegia{| \Psi_n(\hat{z}_n) | +{}} \frac{3}{\lambda} |\zeta_n(\sigma)|_X^2 + \frac{1}{\lambda} | \xi(\sigma) |_X^2 +\frac{1}{\lambda} |\hat{z}_n|_X^2 + \frac{\lambda}{2} c_0^2 + d_0
\\[2ex]
                      & \leq & \ds\frac{16}{\lambda}|\xi(\sigma)|_X^2 + \hat{M}_2, \mbox{ for $ \mu $-a.e. $ \sigma \in S $ and any $ n \in \N $,}
\end{array}
\end{equation}
where
\begin{equation*}
\hat{M}_2 := \frac{3}{\lambda} \hat{M}_1 + \frac{\lambda}{2} c_0^2 + d_0 + \sup_{n \in \N} \left( \Psi_n(\hat{z}_n) + \frac{1}{\lambda} |\hat{z}_n|_X^2 \right).
\end{equation*}
\end{comment}
%%%%%%%%%%%%%%%%%%
In the light of \eqref{axLem21}--\eqref{axLem22}, we can apply Lebesgue's dominated convergence theorem, to derive that:
\begin{equation}\label{axLem23}
\hat{\Psi}_n(\zeta_n) \to \hat{\Psi}(\zeta) \mbox{ as $ n \to \infty $.}
\end{equation}
\piegia{Then,} \eqref{axLem16}--\eqref{axLem18}, \eqref{axLem23} and the previous (ii-1) imply the validity \piegia{of~(ii-2)}. \hfill $ \Box $
\medskip

Now the Key-Lemmas 1--3 are proved as follows.
\\
\textbf{Proof of Key-Lemma \ref{keyLem0}.}\ \ First, we show that $ \partial \Phi_{\varepsilon,\delta}^{\lambda} \subset \mathcal{A}_{\varepsilon,\delta}^{\lambda} $ in $ \mathscr{H} \times \mathscr{H} $. Let us assume that $ [u,u_{\Gamma}] \in D(\partial \Phi_{\varepsilon,\delta}^{\lambda}) $ and $ [u^{*},u_{\Gamma}^{*}] \in \partial \Phi_{\varepsilon,\delta}^{\lambda} (u,u_{\Gamma}) $ in $ \mathscr{H} $. Besides, let us take arbitrary $ \tau > 0 $ and $ [z,z_{\Gamma}] \in \mathscr{V}_{\varepsilon} $,  to compute that:
\begin{equation*}
\begin{array}{lll}
& \multicolumn{2}{l}{(u^*, z)_{L^2(\Omega)} + (u_{\Gamma}^*, z_{\Gamma})_{L^2(\Gamma)}}
\\[2ex]
\le & \multicolumn{2}{l}{\ds\frac{1}{\tau} \left\{ \Phi_{\varepsilon,\delta}^{\lambda} (u + \tau z, u_{\Gamma} + \tau z_{\Gamma}) - \Phi_{\varepsilon,\delta}^{\lambda} (u,u_{\Gamma})\right\}}
\\[2ex]
= & \multicolumn{2}{l}{\ds\frac{1}{\tau} \int_{\Omega} \left( f_{\delta} (\nabla (u + \tau z)) + \frac{\kappa^2}{2} |\nabla (u + \tau z) |^2 - f_{\delta}(\nabla u) - \frac{\kappa^2}{2} |\nabla u|^2 \right) \, dx} 
\\[2ex]
& & \ds + \frac{1}{2\tau} \int_{\Gamma} \left( |\nabla_{\Gamma} (\varepsilon (u_{\Gamma} + \tau z_{\Gamma})) |^2 - | \nabla_{\Gamma}  (\varepsilon u_{\Gamma}) |^2 \right) \, d\Gamma
\\[2ex]
& & \ds + \frac{1}{\tau} \int_{\Omega} \left( B^{\lambda}( u + \tau z) - B^{\lambda}(u) \right) \, dx + \frac{1}{\tau} \int_{\Gamma} \left( B_{\Gamma}^{\lambda} ( u_{\Gamma} + \tau z_{\Gamma} ) - B_{\Gamma}^{\lambda} (u_{\Gamma}) \right) \, d\Gamma
\\[2ex]
\noalign{\allowbreak}
\to & \multicolumn{2}{l}{\ds\int_{\Omega} \left( \nabla f_{\delta}(\nabla u) + \piegia{\kappa^2} \, \nabla u \right) \cdot \nabla z \, dx + \int_{\Gamma} \nabla_{\Gamma} (\varepsilon u_{\Gamma}) \cdot \nabla_{\Gamma} (\varepsilon z_{\Gamma}) \, d\Gamma}
\\[2ex]
& & \ds+ \int_{\Omega} \beta^{\lambda}(u) z \, dx + \int_{\Gamma} \beta_{\Gamma}^{\lambda} (u_{\Gamma}) z_{\Gamma} \, d\Gamma \mbox{ as $ \tau \downarrow 0 $. }
\end{array}
\end{equation*}
Since the choice of $ [z,z_{\Gamma}] \in \mathscr{V}_{\varepsilon} $ is arbitrary, we have
\begin{equation}\label{key0.01}
\begin{array}{lll}
& \multicolumn{2}{l}{\ds \left( u^* - \beta^{\lambda}(u), z \right)_{L^2(\Omega)} + \left( u_{\Gamma}^{*} - \beta_{\Gamma}^{\lambda}(u_{\Gamma}) , z_{\Gamma} \right)_{L^2(\Gamma)}}
\\[2ex]
= & \multicolumn{2}{l}{\ds \int_{\Omega} \left( \nabla f_{\delta} (\nabla u) + \kappa^2 \nabla u \right) \cdot \nabla z \, dx + \int_{\Gamma} \nabla_{\Gamma} (\varepsilon u_{\Gamma}) \cdot \nabla_{\Gamma} (\varepsilon z_{\Gamma}) \, d\Gamma,}
\\[2ex]
\multicolumn{3}{c}{\mbox{for any $ [z,z_{\Gamma}] \in \mathscr{V}_{\varepsilon} $.}}
\end{array}
\end{equation}
Here, taking any $ \varphi_0 \in H_0^1(\Omega) $ and putting $ [z,z_{\Gamma}] = [\varphi_0, 0] \in \mathscr{V}_{\varepsilon} $ in \eqref{key0.01}, 
\begin{equation*}
\begin{array}{c}
\ds(u^* - \beta^{\lambda}(u), \varphi_0)_{L^2(\Omega)} = \int_{\Omega} \bigl( \nabla f_{\delta} (\nabla u) + \kappa^2 \nabla u \bigr) \cdot \nabla \varphi_0 \, dx,
\\[2ex] 
\mbox{for any $ \varphi_0 \in H_0^1(\Omega) $,}
\end{array}
\end{equation*}
which implies
\begin{equation}\label{key0.02}
-\mathrm{div} \left( \nabla f_{\delta} (\nabla u) + \kappa^2 \nabla u \right) = u^* - \beta^{\lambda}(u) \in L^2(\Omega) \mbox{ in $ \mathscr{D}'(\Omega) $.}
\end{equation}
Additionally, with Remark \ref{Rem.sOps} (Fact\,1)--(Fact\,3) and \eqref{key0.01}--\eqref{key0.02} in mind, we can see that:
\begin{equation*}
\begin{array}{lll}
& \multicolumn{2}{l}{(u_{\Gamma}^* - \beta_{\Gamma}^{\lambda}(u_{\Gamma}), z_{\Gamma})_{L^2(\Gamma)}} 
\\[2ex]
= & \multicolumn{2}{l}{\ds\int_{\Omega} \left( \nabla f_{\delta}(\nabla u) + \kappa^2 \nabla u \right) \cdot \nabla z \, dx - \left( u^* - \beta^{\lambda}(u), z \right)_{L^2(\Omega)} + \int_{\Gamma} \nabla_{\Gamma} (\varepsilon u_{\Gamma}) \cdot \nabla_{\Gamma} (\varepsilon z_{\Gamma}) \, d\Gamma}
\\[2ex]
= & \multicolumn{2}{l}{\ds{}_{H^{- \frac{1}{2}}(\Gamma)} \left\langle \left[\big( \nabla f_{\delta} (\nabla u) + \kappa^2 \nabla u \big) \cdot n_{\Gamma} \right]_{_{\Gamma}}, z_{\Gamma} \right\rangle_{H^{\frac{1}{2}}(\Gamma)} + {}_{H^{- 1}(\Gamma)} \left\langle - {\mit \Delta}_{\Gamma} (\varepsilon u_{\Gamma}), \varepsilon z_{\Gamma} \right\rangle_{H^1(\Gamma)},}
\\[2ex]
\multicolumn{3}{c}{\mbox{for any $ [z,z_{\Gamma}] \in \mathscr{V}_{\varepsilon} $.}}
\end{array}
\end{equation*}
This identity leads to:
\begin{equation}\label{key0.03}
- {\mit \Delta}_{\Gamma} (\varepsilon^2 u_{\Gamma}) + \left[ \bigl( \nabla f_{\delta} (\nabla u) + \kappa^2 \nabla u \bigr) \cdot n_{\Gamma} \right]_{_{\Gamma}} = u_{\Gamma}^* - \beta_{\Gamma}^{\lambda}(u_{\Gamma}) \in L^2(\Gamma) \mbox{ in $ H^{-1}(\Gamma) $. }
\end{equation}
As a consequence of \eqref{key0.01}--\eqref{key0.03}, we obtain that:
\begin{equation}\label{key0.04}
[u,u_{\Gamma}] \in \mathscr{D}_{\varepsilon,\delta}^{\lambda} \mbox{ and } [u^*,u_{\Gamma}^*] \in \mathcal{A}_{\varepsilon,\delta}^{\lambda} [u,u_{\Gamma}] \mbox{ in $ \mathscr{H} $.}
\end{equation}

Secondly, we show that $ \mathcal{A}_{\varepsilon,\delta}^{\lambda} \subset \partial \Phi_{\varepsilon,\delta}^{\lambda} $ in $ \mathscr{H} \times \mathscr{H} $. Let us assume that $ [u,u_{\Gamma}] \in \mathscr{D}_{\varepsilon,\delta}^{\lambda} $ and $ [u^*,u_{\Gamma}^{*}] \in \mathcal{A}_{\varepsilon,\delta}^{\lambda} [u,u_{\Gamma}] $ in $ \mathscr{H} $, and let us take an arbitrary $ [z,z_{\Gamma}] \in \mathscr{V}_{\varepsilon} $. Then, taking into account Remark \ref{Rem.sOps} (Fact\,1)--(Fact\,3)\piegia{,} (A4)--(A5) \piegia{and the convexity of the squared norm}, we compute that:
\begin{equation*}
\begin{array}{lll}
& \multicolumn{2}{l}{([u^*,u_{\Gamma}^*], [z,z_{\Gamma}] - [u,u_{\Gamma}])_{\mathscr{H}}}
\\[2ex]
\noalign{\allowbreak}
= & \multicolumn{2}{l}{\ds\left( -\mathrm{div}\left( \nabla f_{\delta}(\nabla u) + \kappa^2 \nabla u \right) + \beta^{\lambda}(u), z - u\right)_{L^2(\Omega)}}
\\[2ex]
& & \ds + \left( - {\mit \Delta}_{\Gamma} (\varepsilon^2 u_{\Gamma}) + \left[ \left( \nabla f_{\delta} (\nabla u) + \kappa^2 \nabla u \right) \cdot n_{\Gamma} \right]_{_{\Gamma}} + \beta_{\Gamma}^{\lambda}(u_{\Gamma}), z_{\Gamma} - u_{\Gamma} \right)_{L^2(\Gamma)}
\\[2ex]
\noalign{\allowbreak}
= & \multicolumn{2}{l}{\ds\int_{\Omega} \left( \nabla f_{\delta}(\nabla u) + \kappa^2 \nabla u \right) \cdot \nabla (z - u) \, dx + \int_{\Omega} \beta^{\lambda}(u)(z - u) \, dx}
\\[2ex]
& & + \ds\int_{\Gamma} \nabla_{\Gamma} (\varepsilon u_{\Gamma} ) \cdot \nabla_{\Gamma} \left( \varepsilon (z_{\Gamma} - u_{\Gamma}) \right) \, d\Gamma + \int_{\Gamma} \beta_{\Gamma}^{\lambda}(u_{\Gamma})(z_{\Gamma} - u_{\Gamma}) \, d\Gamma
\\[2ex]
\leq & \multicolumn{2}{l}{\ds\int_{\Omega} \left( f_{\delta}(\nabla z) \piegia{{}- f_{\delta}(\nabla u) + \frac{\kappa^2}{2} |\nabla z|^2{}} - \frac{\kappa^2}{2} |\nabla u|^2 \right) \, dx + \int_{\Omega} \left( B^{\lambda}(z) - B^{\lambda}(u) \right) \, dx}
\\[2ex]
& & \ds + \frac{1}{2} \int_{\Gamma} \left( |\nabla_{\Gamma}(\varepsilon z_{\Gamma})|^2 - |\nabla_{\Gamma}(\varepsilon u_{\Gamma})|^2\right) \, d\Gamma + \int_{\Gamma} \left( B_{\Gamma}^{\lambda}(z_{\Gamma}) - B_{\Gamma}^{\lambda}(u_{\Gamma}) \right) \, d\Gamma
\\[2ex]
\leq & \multicolumn{2}{l}{\ds\Phi_{\varepsilon,\delta}^{\lambda}(z,z_{\Gamma}) - \Phi_{\varepsilon,\delta}^{\lambda}(u,u_{\Gamma}) \mbox{, for any $ [z,z_{\Gamma}] \in \mathscr{V}_{\varepsilon} $.}}
\end{array}
\end{equation*}
It implies that:
\begin{equation}\label{key0.05}
[u,u_{\Gamma}] \in D(\partial \Phi_{\varepsilon,\delta}^{\lambda}) \mbox{ and } [u^*,u_{\Gamma}^*] \in \partial \Phi_{\varepsilon,\delta}^{\lambda}(u,u_{\Gamma}) \mbox{ in $ \mathscr{H} $.} 
\end{equation}

Thus, we conclude this lemma by \eqref{key0.04}--\eqref{key0.05}. \hfill $ \Box $
\medskip

\noindent
\textbf{Proof of Key-Lemma \ref{keyLem2}.}\ \ First, we verify the condition of lower-bound. Let $ [\check{u}, \check{u}_{\Gamma}] \in \mathscr{H} $ and $ \left\{ [\check{u}_{n}, \check{u}_{\Gamma, n}] \right\}_{n = 1}^{\infty} \subset \mathscr{H} $ be such that:
\begin{equation}\label{key2.1}
[\check{u}_{n}, \check{u}_{\Gamma, n}] \to [\check{u}, \check{u}_{\Gamma}] \mbox{ weakly in $ \mathscr{H} $ as $ n \to \infty $.}
\end{equation}
Then, we may suppose $ \ds{\varliminf_{n \to \infty} \Phi_{\varepsilon_n, \delta_n}^{\lambda_n}(\check{u}_{n}, \check{u}_{\Gamma, n})} < \infty $, 
because \piegia{the other} case is trivial. \piegia{Hence}, 
there exists a subsequence $ \{ n_k \}_{k = 1}^{\infty} \subset \{ n \}_{n = 1}^{\infty} $, such that:
\begin{equation}\label{key2.2}
\check{\Phi}_* := \varliminf_{n \to \infty } \Phi_{\varepsilon_n, \delta_n}^{\lambda_n} (\check{u}_{n}, 
\check{u}_{\Gamma, n}) = \lim_{k \to \infty} \Phi_{\varepsilon_{n_k}, \delta_{n_k}}^{\lambda_{n_k}} (\check{u}_{n_k}, \check{u}_{\Gamma, n_k}) < \infty.
\end{equation}
Here, from \eqref{ap.Phi}, it can be seen that $ \left\{ [\check{u}_{n_k}, \check{u}_{\Gamma, n_k}] \right\}_{k = 1}^{\infty} $ 
is bounded in $ \mathscr{V}_0 $ (resp.\ $ \mathscr{V}_{\varepsilon_0} $), 
if $ \varepsilon_0 = 0 $ (resp.\ if $ \varepsilon_0 > 0 $). 
So, by invoking (Fact\,6) and \eqref{key2.1}, \piegia{and taking} more subsequences if necessary, we can further suppose that:
\begin{equation}\label{key2.3}
\left\{
\begin{array}{l}
[\check{u}_{n_k}, \check{u}_{\Gamma, n_k}] \to [\check{u}, \check{u}_{\Gamma}] \mbox{ in $ \mathscr{H} $ and weakly in $ \mathscr{V}_0$,}
\\[0.5ex]
\check{u}_{n_k}(x) \to \check{u}(x) \mbox{ a.e. $ x \in \Omega $,}
%\\[0.5ex]
%\check{u}_{\Gamma, n_k}(x_{\Gamma}) = (\check{u}_{n_k})_{|_{\Gamma}}(x_{\Gamma}) \to \check{u}_{|_{\Gamma}}(x_{\Gamma}) = \check{u}_{\Gamma}(x_{\Gamma}) \mbox{ a.e. $ x_{\Gamma} \in \Gamma $,}
\end{array}
\right. \mbox{ as $ k \to \infty $,}
\end{equation}
and in particular, if $ \varepsilon_0 > 0 $, then
\begin{equation}\label{key2.4}
\check{u}_{\Gamma, n_k} \to \check{u}_{\Gamma} \mbox{ weakly in $ H^1(\Gamma) $ as $ k \to \infty $.}
\end{equation}
Additionally, from (A4) and Lemma \ref{auxLem}, we can infer that:
\begin{description}
\item[(Fact\,8)]The sequence of convex functions:
\begin{equation*}
\left\{ \omega \in L^2(\Omega)^N \mapsto \int_{\Omega} f_{\delta_n}(\omega) \, dx \in [0,\infty] \right\}_{n=1}^{\infty}
\end{equation*}
converges to the convex function of $ L^1 $-norm:
\begin{equation*}
\omega \in L^2(\Omega)^N \mapsto \int_{\Omega} |{}\omega{}| \, dx \in [0,\infty),
\end{equation*}
on the Hilbert space $ L^2(\Omega)^N $, in the sense of \piegia{Mosco}, as $ n \to \infty $.
\end{description}
In the light of Remark \ref{MY}, \eqref{key2.1}--\eqref{key2.4}, Fatou's lemma and the above (Fact\,8), the condition of lower-bound is verified as follows:
\begin{equation*}
\begin{array}{lll}
\check{\Phi}_* & \multicolumn{2}{l}{\geq \ds\varliminf_{k \to \infty} \int_{\Omega} f_{\delta_{n_k}}(\nabla \check{u}_{n_k}) \, dx + \frac{\kappa^2}{2} \varliminf_{k \to \infty} \int_{\Omega} |\nabla \check{u}_{n_k}|^2 \, dx + \varliminf_{k \to \infty} \int_{\Omega} B^{\lambda_{n_k}}(\check{u}_{n_k}) \, dx}
\\[2ex]
&  & + \ds\frac{1}{2} \varliminf_{k \to \infty} \int_{\Gamma} \left|\nabla_{\Gamma} (\varepsilon_{n_k} \check{u}_{\Gamma, n_k}) \right|^2 \, d\Gamma + \varliminf_{k \to \infty} \int_{\Gamma} B_{\Gamma}^{\lambda_{n_k}}(\check{u}_{\Gamma,n_k}) \, d\Gamma
\\[2ex]
& \multicolumn{2}{l}{\geq \Phi_{\varepsilon_0}(\check{u},\check{u}_{\Gamma}).}
\end{array}
\end{equation*}

Next, we verify the optimality condition. Let us fix any $ [\hat{u}_0, \hat{u}_{\Gamma, 0}] \in D(\Phi_{\varepsilon_0}) $, and let us take a sequence $ \{ \omega_i \}_{ i = 1 }^{ \infty } \subset H^1(\Omega) $ in the following way:
\begin{equation}\label{key2.5}
\left\{ \hspace{-2.5ex}\parbox{12.5cm}{
	\vspace{-2ex}
	\begin{itemize}
	\item if $ \varepsilon_0 > 0 $, then $ \{ \omega_i \}_{i = 1}^\infty = \{ \hat{u}_0 \} $;
	\item if $ \varepsilon_0 = 0 $, then $ \{ \omega_i \}_{i = 1}^\infty \subset \piegia{{}C^1{}}(\overline{\Omega}) $ is such that $ \omega_i \to \hat{u}_0 $ in $ H^1(\Omega) $, and pointwisely a.e. in $ \Omega $, as $ i \to \infty $.
	\vspace{-2ex}
	\end{itemize}
} \right.
\end{equation}
%\begin{equation}\label{key2.5}
%\left\{\begin{array}{l}
%\omega_i \to \hat{u}_0 \mbox{ \ in $ H^1(\Omega) $, and pointwisely a.e. in $ \Omega $, as $ i \to \infty $,}
%\\[0.5ex]
%(\omega_{i})_{|_{\Gamma}} = \hat{u}_{\Gamma,0} (= (\hat{u}_0)_{_{\Gamma}}) \mbox{ in $ H^{\frac{1}{2}}(\Gamma) $, for any $ i \in \N $.}
%\end{array} \right.
%\end{equation}
Besides, let us define a sequence $ \{ \hat{\varphi}_i \}_{ i = 1}^{ \infty } \subset H^1(\Omega) $, by putting:
\begin{equation*}
\hat{\varphi_i} := \mathcal{T}_{B} \omega_i \mbox{ \ in $ H^1(\Omega) $, for $ i = 1,2,3,\dots $\piegia{,}}
\end{equation*}
\piegia{where the projection function $\mathcal{T}_{B}$ is defined by \eqref{T_B^rho}.} 
Then, \piegia{in view of \eqref{key2.5}, and taking a subsequence if necessary, we have that} 
\begin{equation}\label{key2.6}
\left\{
\begin{array}{l}
\hat{\varphi}_i \to \hat{u}_0 \mbox{ in $ H^1(\Omega) $, and pointwisely a.e. in $ \Omega$, as $ i \to \infty $}\\[0.5ex]
(\hat{\varphi}_{i})_{|_\Gamma}\to \hat{u}_{\Gamma, 0} \mbox{ in $ H^{\frac{1}{2}}(\Gamma) $, and pointwisely a.e. on $ \Gamma $, as $ i \to \infty $}\\[0.5ex]
\{ (\hat{\varphi}_{i})_{|_\Gamma} \}_{i = 1}^\infty \subset H^1(\Gamma), \mbox{ and } (\varepsilon_0 \hat{\varphi}_{i})_{|_\Gamma}\to \varepsilon_0 \hat{u}_{\Gamma, 0} \mbox{ in $ H^{1}(\Gamma) $, as $ i \to \infty $.}
\end{array}
\right.
\end{equation}
Also, invoking (A4) and Lebesgue's dominated convergence theorem, we can configure a sequence $ \{ n_i \}_{i = 0}^{\infty} \subset \N $, such that $ 1 =: n_0 < n_1 < n_2 < n_3 < \cdots < n_i \uparrow \infty $, as $ i \to \infty $, and for any $ i \in \N \cup \{ 0 \} $,
\begin{equation}\label{key2.8}
\left\{\begin{array}{l}
\ds\sup_{n \geq n_i} \left| f_{\delta_n}(\nabla \hat{\varphi}_i) - | \nabla \piegia{\hat{\varphi}_i}| \/\right|_{L^1(\Omega)} < \frac{1}{2^{i + 1}},
\\[2ex]
\ds\sup_{n \geq n_i}\left( \left| B^{\lambda_n}(\hat{\varphi}_i) - B(\hat{\varphi}_i) \right|_{L^1(\Omega)} + \left| B_{\Gamma}^{\lambda_n}(\hat{\varphi}_i) - B_{\Gamma}(\hat{\varphi}_i) \right|_{L^1(\Gamma)} \right) < \frac{1}{2^{i}},
\\[2ex]
\ds\left( \sup_{ n \geq n_i } (\varepsilon_n^2 - \varepsilon_0^2) \right) |\nabla \hat{\varphi}_i|_{L^2(\Omega)}^2 < \frac{1}{2^{ i + 1}}.
\end{array} \right.
\end{equation}
Based on these, let us define:
\begin{equation}\label{key2.9}
[\hat{u}_n, \hat{u}_{\Gamma, n}] := [\hat{\varphi}_i, (\hat{\varphi}_i)_{|_{\Gamma}}] \mbox{, if $ n_i \leq n < n_{i + 1} $, for some $ i \in \N \cup \{ 0 \} $.}
\end{equation}
Then, with \piegia{condition (a2) in (A1)}, \eqref{key2.6} and Lebesgue's dominated convergence theorem in mind, one can see that:
\begin{equation}\label{key2.9-1}
	\lim_{n \to \infty} \int_{\Omega} B(\hat{u}_n) \, dx = \int_{\Omega} B(\hat{u}_0) \, dx, \mbox{ and } 	\lim_{n \to \infty} \int_{\Gamma} B_\Gamma(\hat{u}_{\Gamma, n}) \, d\Gamma = \int_{\Gamma} B_\Gamma(\hat{u}_{\Gamma, 0}) \, dx.
\end{equation}
Taking into account \eqref{key2.5}--\eqref{key2.9-1}, we obtain that:
\begin{equation*}
\begin{array}{lll}
& \multicolumn{2}{l}{\ds\left| \Phi_{\varepsilon_n, \delta_n}^{\lambda_n} (\hat{u}_n, \hat{u}_{\Gamma, n}) - \Phi_{\varepsilon_0}(\hat{u}_0, \hat{u}_{\Gamma, 0}) \right|}
\\[2ex]
\leq & \multicolumn{2}{l}{\ds\int_{\Omega} \bigl| f_{\delta_n}(\nabla \hat{u}_n) - |\nabla \hat{u}_n| \bigr| \, dx + \left| \int_{\Omega} \left(|\nabla \hat{u}_n| - |\nabla \hat{u}_0| \right) \, dx \right|}
\\[2ex]
& & + \ds\frac{\kappa^2}{2} \int_{\Omega} \left| |\nabla \hat{u}_n|^2 - |\nabla \hat{u}_0|^2 \right| \, dx + \frac{1}{2} \left| \int_{\Gamma} \left(|\nabla_{\Gamma}(\varepsilon_n \piegia{\hat{u}_{\Gamma, n}})|^2 - | \nabla_{\Gamma} (\varepsilon_0 \hat{u}_{\Gamma, 0}) |^2 \right) \, d\Gamma \right|
\\[2ex]
& & + \ds\left| \int_{\Omega}  B^{\lambda_n}(\hat{u}_n) \, dx - \int_{\Omega} B(\hat{u}_n) \, dx \right| + \left| \int_{\Omega}  B(\hat{u}_n) \, dx - \int_{\Omega} B(\hat{u}_0) \, dx \right|
\\[2ex]
& &+ \ds\left| \int_{\Gamma} B_{\Gamma}^{\lambda_n}(\hat{u}_{\Gamma, n}) \, d\Gamma - \int_{\Gamma} B_{\Gamma}(\hat{u}_{\Gamma, n}) \, d\Gamma \right| + \left| \int_{\Gamma} B_{\Gamma}(\hat{u}_{\Gamma, n}) \, d\Gamma - \int_{\Gamma} B_{\Gamma}(\hat{u}_{\Gamma, 0}) \, d\Gamma \right|
\\[2ex]
\leq & \multicolumn{2}{l}{\ds\int_{\Omega} \bigl| f_{\delta_n}(\nabla \hat{u}_n) - |\nabla \hat{u}_n| \bigr| \, dx + | \nabla (\hat{u}_n - \hat{u}_0) |_{L^1(\Omega)^N}}
\\[2ex]
& & + \ds\frac{\kappa^2}{2} \bigl| |\nabla \hat{u}_n| + |\nabla \hat{u}_0| \bigr|_{L^2(\Omega)} \left| \nabla (\hat{u}_n - \hat{u}_0) \right|_{L^2(\Omega)^N} + \frac{1}{2} |\varepsilon_n^2 - \varepsilon_0^2| \left| \nabla_{\Gamma} \hat{u}_{\Gamma, n} \right|_{L^2(\Gamma)^N}^2
\\[2ex]
& & + \ds\frac{1}{2} \bigl| | \nabla_{\Gamma} ( \varepsilon_0 \hat{u}_{\Gamma, n} ) | + | \nabla_{\Gamma} ( \varepsilon_0 \hat{u}_{\Gamma, 0} ) | \bigr|_{L^2(\Gamma)} \left| \nabla_{\Gamma} (\varepsilon_0 (\hat{u}_{\Gamma, n} - \hat{u}_{\Gamma, 0})) \right|_{L^2(\Gamma)^N}
\\[2ex]
& & + \ds\left| B^{\lambda_n}(\hat{u}_n) - B(\hat{u}_n) \right|_{L^1(\Omega)} + \left| B_{\Gamma}^{\lambda_n}(\hat{u}_{\Gamma, n}) - B_{\Gamma}(\hat{u}_{\Gamma, 0}) \right|_{L^1(\Gamma)}
\\[2ex]
& & + \ds\left| \int_{\Omega} B(\hat{u}_n) \, dx - \int_{\Omega} B(\hat{u}_0) \, dx \right| + \left| \int_{\Gamma} B_{\Gamma}(\hat{u}_{\Gamma, n}) \, d\Gamma - \int_{\Gamma} B_{\Gamma}(\hat{u}_{\Gamma, 0}) \, d\Gamma \right|
%\\[2ex]
\end{array}
\end{equation*}
\begin{equation*}
\begin{array}{lll}
	\leq & \multicolumn{2}{l}{\ds\frac{1}{2^{i-1}} + \left( 1 + \mathcal{L}^{N}(\Omega)^{\frac{1}{2}} + \hat{\Phi}_* \right) \left( \rule{0pt}{18pt} | \hat{u}_n - \hat{u}_0 |_{H^1(\Omega)} + | \varepsilon_0 (\hat{u}_{\Gamma, n} - \hat{u}_{\Gamma, n}) |_{H^1(\Gamma)} \right.}
\\[2ex]
& \multicolumn{2}{l}{\qquad\ds\left.  + \left| \int_{\Omega} B(\hat{u}_n) \, dx - \int_{\Omega} B(\hat{u}_0) \, dx \right| + \left| \int_{\Gamma} B_{\Gamma}(\hat{u}_{\Gamma, n}) \, d\Gamma - \int_{\Gamma} B_{\Gamma}(\hat{u}_{\Gamma, 0}) \, d\Gamma \right| \right),}
\\[2ex]
\multicolumn{3}{c}{\mbox{for any $ i \in \N \cup \{ 0 \} $ and $ n \geq n_i $,}}
\end{array}
\end{equation*}
where
\begin{equation*}
	\hat{\Phi}_* := \sup_{ n \in \N } \left( \frac{\kappa^2}{2} \bigl| | \nabla \hat{u}_n | + | \nabla \hat{u}_0 | \bigr|_{L^2(\Omega)} +\frac{1}{2} \bigl| |\nabla_{\Gamma} (\varepsilon_0 u_{\Gamma,n})| +|\nabla_{\Gamma} (\varepsilon_0 u_{\Gamma,0})| \bigr|_{L^2(\Gamma)}\right).
\end{equation*}
This implies that the sequence $ \{ [\hat{u}_n, \hat{u}_{\Gamma, n}] \} \subset H^1(\Omega) \times H^1(\Gamma) $ is the required sequence to verify the condition of optimality. \hfill$ \Box $
\medskip

By a similar demonstration technique, we also see the following Corollary.
\begin{cor}\label{cor.1}
Let $ \{ \varepsilon_n \}_{n=1}^{\infty} \subset [0,\infty) $ be arbitrary sequence such that $ \varepsilon_{n} \to \varepsilon_0 $ as $ n \to \infty $. Then, for the sequence $ \{ \Phi_{\varepsilon_n} \}_{n = 1}^{\infty} $ of convex functions, it holds that:
\begin{equation*}
\Phi_{\varepsilon_{n}} \to \Phi_{\varepsilon_0} \mbox{ on $ \mathscr{H} $, in the sense of \piegia{Mosco}, as $ n \to \infty $.}
\end{equation*}
\end{cor}
\medskip

\noindent
\textbf{Proof of Key-Lemma \ref{keyLem1}.}\ \ Let us fix any $ \varepsilon \geq 0 $, and let us define a set-valued map $ \mathcal{A}_{\varepsilon} : \mathscr{H} \to 2^{\mathscr{H}} $, by putting:
\begin{equation}\label{key1.01}
D(\mathcal{A}_{\varepsilon}) := \left\{\begin{array}{l|l} [u,u_{\Gamma}] \in \mathscr{V}_{\varepsilon} & \begin{tabular}{l} \mbox{there exists $ \nu_u \in L^{\infty}(\Omega)^N $ and $ [\xi, \xi_{\Gamma}] \in \mathscr{H} $}
\\[1ex]
\mbox{such that \eqref{key1.1}--\eqref{key1.2} \piegia{hold}.} \end{tabular} \end{array}\right\},
\end{equation}
and 
\begin{equation}\label{key1.02}
\begin{array}{rcl}
[u,u_{\Gamma}] & \in & D(\mathcal{A}_{\varepsilon}) \subset \mathscr{H} \mapsto \mathcal{A}_{\varepsilon}[u,u_{\Gamma}]
\\[2ex]
& := &\left\{\begin{array}{l|l} [u^*,u_{\Gamma}^*] \in \mathscr{H} 
& \begin{tabular}{l} 
\eqref{key1.3} holds, for some $ \nu_u \in L^{\infty}(\Omega)^N $ and \\[1ex]
$ [\xi,\xi_{\Gamma}] \in \mathscr{H} $, satisfying \eqref{key1.1}--\eqref{key1.2} \end{tabular} \end{array}\right\}.
\end{array}
\end{equation}
Then, the assertion of Key-Lemma \ref{keyLem1} can be rephrased as follows:
\begin{equation}\label{key1.03}
\partial \Phi_{\varepsilon} = \mathcal{A}_{\varepsilon} \mbox{ in $ \mathscr{H} \times \mathscr{H} $.}
\end{equation}
This coincidence will be obtained as a consequence of the following {\it Claims $ \# $1--$ \# $2}.
\medskip

\noindent
{\it Claim $ \# $1: $ \mathcal{A}_{\varepsilon} \subset \partial \Phi_{\varepsilon} $ in $ \mathscr{H} \times \mathscr{H} $.}

Let us assume that $ [u,u_{\Gamma}] \in D(\mathcal{A}_{\varepsilon}) $ and $ [u^*,u_{\Gamma}^*] \in \mathcal{A}_{\varepsilon}[u,u_{\Gamma}] $ in $ \mathscr{H} $. Then, from (A1), Remark \ref{Rem.sOps} (Fact\,1)--(Fact\,3) and Remark \ref{exConvex}, and \eqref{key1.01}--\eqref{key1.02}, it is inferred that:
\begin{equation*}
\begin{array}{rcl}
& \multicolumn{2}{l}{\bigl([u^*,u_{\Gamma}^*], [z,z_{\Gamma}] - [u,u_{\Gamma}]\bigr)_{\mathscr{H}}}
\\[2ex]
= & \multicolumn{2}{l}{\ds\left(-\mathrm{div}(\nu_u + \kappa^2 \nabla u) + \xi, z - u \right)_{L^2(\Omega)}} 
\\[2ex]
& & + \ds\left( -{\mit \Delta}_{\Gamma} (\varepsilon^2 u_{\Gamma}) + \left[\bigl(\nu_u + \kappa^2 \nabla u \bigr) \cdot n_{\Gamma} \right]_{_{\Gamma}} + \xi_{\Gamma}, z_{\Gamma} - u_{\Gamma} \right)_{L^2(\Gamma)}
\\[2ex]
= & \multicolumn{2}{l}{\ds\int_{\Omega} (\nu_u + \kappa^2 \nabla u) \cdot \nabla(z - u) \, dx + \int_{\Omega} \xi (z - u) \, dx}
\\[2ex]
& & + \ds\int_{\Gamma} \nabla_{\Gamma} (\varepsilon u_{\Gamma}) \cdot \nabla_{\Gamma} ( \varepsilon (z_{\Gamma} - u_{\Gamma})) \, d\Gamma + \int_{\Gamma} \xi_{\Gamma} (z_{\Gamma} - u_{\Gamma}) \, d\Gamma
\\[2ex]
\leq & \multicolumn{2}{l}{\ds\int_{\Omega} \left( | \nabla z | \piegia{{}- | \nabla u | + \frac{\kappa^2}{2} |{}}\nabla z|^2   - \frac{\kappa^2}{2} |\nabla u|^2 \right) \, dx + \int_{\Omega} \left( B(z) - B(u) \right) \, dx}
\\[2ex]
& & + \ds\frac{1}{2} \int_{\Gamma} \left( |\nabla_{\Gamma}(\varepsilon z_{\Gamma})|^2 - |\nabla_{\Gamma}(\varepsilon u_{\Gamma})|^2\right) \, d\Gamma + \int_{\Gamma} \left( B_{\Gamma}(z_{\Gamma}) - B_{\Gamma}(u_{\Gamma}) \right) \, d\Gamma
\\[2ex]
\leq & \multicolumn{2}{l}{\ds\Phi_{\varepsilon} (z,z_{\Gamma}) - \Phi_{\varepsilon} (u,u_{\Gamma}) \mbox{, for any $ [z,z_{\Gamma}] \in \mathscr{V}_{\varepsilon} $.}}
\end{array}
\end{equation*}
Thus, we have:
\begin{equation*}
[u,u_{\Gamma}] \in D(\partial \Phi_{\varepsilon}) \mbox{ and } [u^*,u_{\Gamma}^*] \in \partial \Phi_{\varepsilon}(u,u_{\Gamma}) \mbox{ in $ \mathscr{H} $},
\end{equation*}
and we can say that:
\begin{equation*}
\mathcal{A}_{\varepsilon} \subset \partial \Phi_{\varepsilon}(u,u_{\Gamma}) \mbox{ in $ \mathscr{H} \times \mathscr{H} $.}
\end{equation*}
\medskip

\noindent
{\it Claim $ \# $2: $ (\mathcal{A}_{\varepsilon} + \mathcal{I}_{\mathscr{H}}) \mathscr{H} = \mathscr{H} $.}

Since, $ (\mathcal{A}_{\varepsilon} + \mathcal{I}_{\mathscr{H}}) \mathscr{H} \subset \mathscr{H} $ is trivial, it is sufficient to prove the converse inclusion. Let us take any $ [w,w_{\Gamma}] \in \mathscr{H} $. Then, by applying Minty's theorem, Key-Lemma \ref{keyLem0} and Remark \ref{Rem.sOps} (Fact\,1)--(Fact\,3), we can configure a class of functions 
$ \{ [u_{\delta}^{\lambda}, u_{\Gamma,\delta}^{\lambda}] \,|\, 0 < \delta, \lambda \le 1 \} \subset \mathscr{V}_{\varepsilon} $, by setting:
\begin{equation*}
\left\{ [u_{\delta}^{\lambda}, u_{\Gamma,\delta}^{\lambda}]
 := (\mathcal{A}_{\varepsilon,\delta}^{\lambda} + \mathcal{I}_{\mathscr{H}})^{-1}[w, w_{\Gamma}], ~ 0
 < \delta, \lambda \le 1 \right\} \mbox{ in $ \mathscr{H} $, }
\end{equation*}
i.e.
\begin{equation}\label{key1.04}
[w - u_{\delta}^{\lambda}, w_{\Gamma} - u_{\Gamma,\delta}^{\lambda}]
 = \partial \Phi_{\varepsilon,\delta}^{\lambda}(u_{\delta}^{\lambda}, u_{\Gamma,\delta}^{\lambda})
 \mbox{ in $ \mathscr{H} $, for any $ 0 < \delta, \lambda \le 1 $,}
\end{equation}
and we can see that:
\begin{equation}\label{key1.05}
\begin{array}{lll}
& \multicolumn{2}{l}{\ds\int_{\Omega} \bigl( \nabla f_{\delta}(\nabla u_{\delta}^{\lambda}) + \kappa^2 \nabla u_{\delta}^{\lambda} \bigr) \cdot \nabla z \, dx +  \int_{\Gamma} \nabla_{\Gamma} (\varepsilon u_{\Gamma,\delta}^{\lambda}) \cdot \nabla_{\Gamma} (\varepsilon z_{\Gamma}) \, d\Gamma}
\\[2ex]
& & + \ds\int_{\Omega} \beta^{\lambda}(u_{\delta}^{\lambda}) z \, dx + \int_{\Gamma} \beta_{\Gamma}^{\lambda}(u_{\Gamma,\delta}^{\lambda}) z_{\Gamma} \, d\Gamma
\\[2ex]
= & \multicolumn{2}{l}{\ds\int_{\Omega} ( w - u_{\delta}^{\lambda} ) z \, dx + \int_{\Gamma} ( w_{\Gamma} - u_{\Gamma,\delta}^{\lambda} ) z_{\Gamma} \, d\Gamma,}
\\[2ex]
\multicolumn{3}{c}{\mbox{for any $ [z,z_{\Gamma}] \in \mathscr{V}_{\varepsilon} $ and any $ 0 < \delta, \lambda \le 1 $.}}
\end{array}
\end{equation}
In the variational form \eqref{key1.05}, let us first put $ [z,z_{\Gamma}] = [u_{\delta}^{\lambda},u_{\Gamma,\delta}^{\lambda}] \in \mathscr{V}_{\varepsilon} $ in \eqref{key1.05}. Then, with (A1), (A4)--(A5) and Young's inequality in mind, we deduce that:
\begin{equation}\label{key1.06}
\begin{array}{ll}
\multicolumn{2}{l}{{\ds  \frac{1}{2} \bigl| [u_{\delta}^{\lambda},u_{\Gamma, \delta}^{\lambda}] \bigr|_{\mathscr{H}} + \kappa^2 \left|\nabla u_{\delta}^{\lambda} \right|_{L^2(\Omega)^N}^2 + \left| \nabla_{\Gamma} (\varepsilon u_{\Gamma,\delta}^{\lambda}) \right|_{L^2(\Gamma)^N}^2 \le \frac{1}{2} \bigl| [w,w_{\Gamma}] \bigr|_{\mathscr{H}}^2},}
\\[2ex]
& \multicolumn{1}{c}{\mbox{for any $ 0 < \delta, \lambda \le 1$.}}
\end{array}
\end{equation}

Next, let us \piegia{take} $ [z,z_{\Gamma}] = [\beta^{\lambda}(u_{\delta}^{\lambda}),\beta^{\lambda}(u_{\Gamma,\delta}^{\lambda})] $ in $ \mathscr{V}_{\varepsilon} $. Then, by applying (A1), (A4)--(A5), \cite[Theorem 3.99]{AFP}, \cite[Lemma 4.4]{CC13} and Schwarz's inequality, 
\begin{equation*}
\begin{array}{lll}
& \multicolumn{2}{l}{\ds\left| \beta^{\lambda} (u_{\delta}^{\lambda}) \right|_{L^2(\Omega)}^2 \leq \left(w, \beta^{\lambda}(u_{\delta}^{\lambda})\right)_{L^2(\Omega)} + \left(w_{\Gamma}, \piegia{\beta^{\lambda}}(u_{\Gamma,\delta}^{\lambda})\right)_{L^2(\Gamma)}}
\\[2ex]
\leq & \multicolumn{2}{l}{\ds| w |_{L^2(\Omega)}^2 + \frac{1}{4} \left|\beta^{\lambda}(u_{\delta}^{\lambda}) \right|_{L^2(\Omega)}^2 + a_1^2 | w_{\Gamma} |_{L^2(\Gamma)}^2 + \frac{1}{4a_1^2} \left| \beta^{\lambda}(u_{\Gamma,\delta}^{\lambda}) \right|_{L^2(\Gamma)}^2}
\\[2ex]
\leq & \multicolumn{2}{l}{\ds\left(| w |_{L^2(\Omega)}^2 + a_1^2 | w_{\Gamma} |_{L^2(\Gamma)}^2 \right) + \frac{1}{4} \left| \beta^{\lambda}(u_{\delta}^{\lambda}) \right|_{L^2(\Omega)}^2 + \frac{1}{4a_1^2} \bigl| a_1 | \beta_{\Gamma}^{\lambda}(u_{\Gamma, \delta}^{\lambda}) | + b_1 \bigr|_{L^2(\Gamma)}^2}
\\[2ex]
\leq & \multicolumn{2}{l}{\ds(1 + a_1^2) \bigl| [w,w_{\Gamma}] \bigr|_{\mathscr{H}}^2 + \frac{1}{4} \left| \beta^{\lambda}(u_{\delta}^{\lambda}) \right|_{L^2(\Omega)}^2 + \frac{1}{2} \left| \beta_{\Gamma}^{\lambda}(u_{\Gamma,\delta}^{\lambda}) \right|_{L^2(\Gamma)}^2 + \frac{b_1^2}{2a_1^2} \mathcal{H}^{N - 1}(\Gamma),}
\end{array}
\end{equation*}
so that
\begin{equation}\label{key1.07}
\begin{array}{c}
\ds \frac{3}{4} \left| \beta^{\lambda}(u_{\delta}^{\lambda}) \right|_{L^2(\Omega)}^2 - \frac{1}{2} \left| \beta_{\Gamma}^{\lambda}(u_{\Gamma,\delta}^{\lambda}) \right|_{L^2(\Gamma)}^2 \leq (1 + a_1^2) \bigl| [w,w_{\Gamma}] \bigr|_{\mathscr{H}}^2 + \frac{b_1^2}{2a_1^2} \mathcal{H}^{N - 1}(\Gamma),
\\[2ex]
\mbox{for any $ 0 < \delta, \/ \lambda \leq 1 $.}
\end{array}
\end{equation}
Similarly, putting $ [z, z_{\Gamma}] = [\beta_{\Gamma}^{\lambda}(u_{\delta}^{\lambda}), \beta_{\Gamma}^{\lambda}(u_{\Gamma,\delta}^{\lambda})] \in \mathscr{V}_{\varepsilon} $ in \eqref{key1.05}, and applying (A1), (A4)--(A5), \cite[Theorem 3.99]{AFP}, \cite[Lemma 4.4]{CC13} and Schwarz's inequality, 
\begin{equation*}
\begin{array}{lll}
& \multicolumn{2}{l}{\ds\left| \beta_{\Gamma}^{\lambda} (u_{\Gamma,\delta}^{\lambda}) \right|_{L^2(\Gamma)}^2 \leq \left( w, \beta_{\Gamma}^{\lambda}(u_{\delta}^{\lambda}) \right)_{L^2(\Omega)} + \left(w_{\Gamma}, \beta_{\Gamma}^{\lambda}(u_{\Gamma,\delta}^{\lambda})\right)_{L^2(\Gamma)}}
\\[2ex]
\leq & \multicolumn{2}{l}{\ds\frac{1}{a_0^2} | w |_{L^2(\Omega)}^2 + \frac{a_0^2}{4} \left|\beta_{\Gamma}^{\lambda}(u_{\delta}^{\lambda}) \right|_{L^2(\Omega)}^2 + | w_{\Gamma} |_{L^2(\Gamma)}^2 + \frac{1}{4} \left|\beta_{\Gamma}^{\lambda}(u_{\Gamma,\delta}^{\lambda}) \right|_{L^2(\Gamma)}^2}
\\[2.5ex]
\leq & \multicolumn{2}{l}{\ds\left(\frac{1}{a_0^2} | w |_{L^2(\Omega)}^2 + | w_{\Gamma} |_{L^2(\Gamma)}^2 \right) + \frac{a_0^2}{4} \left| \frac{1}{a_0} \left| \beta^{\lambda}(u_{\delta}^{\lambda}) \right| + \frac{b_0}{a_0} \right|_{L^2(\Omega)}^2 + \frac{1}{4} \left| \beta_{\Gamma}^{\lambda} (u_{\Gamma,\delta}^{\lambda}) \right|_{L^2(\Gamma)}^2}
\\[2.5ex]
\leq & \multicolumn{2}{l}{\ds\left(1 + \frac{1}{a_0^2}\right)\bigl| [w,w_{\Gamma}] \bigr|_{\mathscr{H}}^2 + \frac{1}{2} \left| \beta^{\lambda}(u_{\delta}^{\lambda}) \right|_{L^2(\Omega)}^2 + \frac{1}{4} \left|\beta_{\Gamma}^{\lambda}(u_{\Gamma,\delta}^{\lambda}) \right|_{L^2(\Gamma)}^2 + \frac{b_0^2}{2} \mathcal{L}^{N}(\Omega),}
\end{array}
\end{equation*}
so that
\begin{equation}\label{key1.08}
\begin{array}{c}
\ds - \frac{1}{2} \left| \beta^{\lambda}(u_{\delta}^{\lambda}) \right|_{L^2(\Omega)}^2 + \frac{3}{4} \left| \beta_{\Gamma}^{\lambda}(u_{\Gamma,\delta}^{\lambda}) \right|_{L^2(\Gamma)}^2 \leq \left(1 + \frac{1}{a_0^2} \right)\bigl| [w,w_{\Gamma}] \bigr|_{\mathscr{H}}^2 + \frac{b_0^2}{2}\mathcal{L}^{N}(\Omega),
\\[2ex]
\mbox{for any $ 0 < \delta, \/ \lambda \leq 1 $.}
\end{array}
\end{equation}
Taking the sum of \eqref{key1.07}--\eqref{key1.08}, it follows that:
\begin{equation}\label{key1.09}
\begin{array}{c}
\hspace{-2ex}\ds \frac{1}{4} \left| \left[ \beta^{\lambda}(u_{\delta}^{\lambda}),\beta_{\Gamma}^{\lambda}(u_{\Gamma,\delta}^{\lambda}) \right] \right|_{\mathscr{H}}^2 \leq \left( 2 + \frac{1}{a_0^2} + a_1^2 \right) \bigl| [w,w_{\Gamma}] \bigr|_{\mathscr{H}}^2 + \frac{b_0^2}{2} \mathcal{L}^{N}(\Omega) + \frac{b_1^2}{2 a_1^2} \mathcal{H}^{N - 1}(\Gamma),
\\[2ex]
\mbox{for any $ 0 < \delta, \/ \lambda \leq 1 $.}
\end{array}
\end{equation}

On account of \eqref{key1.06} and \eqref{key1.09}, we find pairs of functions $ [u,u_{\Gamma}] \in \mathscr{V}_{\varepsilon} $ and $ [\xi,\xi_{\Gamma}] \in \mathscr{H} $ and sequences
\begin{equation*}
\left\{\begin{array}{lcl}
\delta_1 > \delta_2 > \delta_3 > \cdots > \delta_n & \downarrow & 0, 
\\[0.5ex]
\lambda_1 > \lambda_2 > \lambda_3 > \cdots > \lambda_n & \downarrow & 0,
\end{array} \right.
\mbox{ as $ n \to \infty $,}
\end{equation*}
such that:
\begin{equation}\label{key1.10}
\left\{\begin{array}{l}
[u_n,u_{\Gamma,n}] := [u_{\delta_n}^{\lambda_n},u_{\Gamma,\delta_n}^{\lambda_n}] \to [u, u_{\Gamma}] \mbox{ in $ \mathscr{H} $ and weakly in $ \mathscr{V}_{\varepsilon} $,}
\\[0.5ex]
\left[\beta^{\lambda_n}(u_n),\beta_{\Gamma}^{\lambda_n}(u_{\Gamma,n}) \right] \to [\xi, \xi_{\Gamma}] \mbox{ weakly in $ \mathscr{H} $,}
\end{array} \right.
\mbox{ as $ n \to \infty $.}
\end{equation}
Here, in the light of Key-Lemma \ref{keyLem2}, \eqref{key1.04} and \eqref{key1.10}, we can apply Remark \ref{Rem.MG} (Fact\,5) to see that:
\begin{equation*}
[w - u,w_{\Gamma} - u_{\Gamma}] \in \partial \Phi_{\varepsilon}(u,u_{\Gamma}) \mbox{ in $ \mathscr{H} $,}
\end{equation*}
and
\begin{equation}\label{key1.11}
\Phi_{\varepsilon,\delta_n}^{\lambda_n}(u_n,u_{\Gamma,n}) \to \Phi_{\varepsilon}(u,u_{\Gamma}) \mbox{ as $ n \to \infty $. }
\end{equation}
Also, by (A1), (A5), Remark \ref{MY} and Remark \ref{Rem.MG} (Fact\,5), we see that:
\begin{equation}\label{key1.10-1}
\xi \in \beta(u) \mbox{ a.e. in $ \Omega $, and } \xi_{\Gamma} \in \beta_{\Gamma}(u_{\Gamma}) \mbox{ a.e. on $ \Gamma $.}
\end{equation}
By virtue of \eqref{key1.10}--\eqref{key1.11}, (A4), Remark \ref{MY}, Lemma \ref{auxLem}, and (Fact\,8), we further compute that:
\begin{equation}\label{key1.12-1}
\begin{array}{lll}
& \multicolumn{2}{l}{\ds\frac{\kappa^2}{2} \int_{\Omega} |\nabla u|^2 \, dx \leq \frac{\kappa^2}{2} \varliminf_{n \to \infty} \int_{\Omega} |\nabla u_n|^2 \, dx \leq \frac{\kappa^2}{2} \varlimsup_{n \to \infty} \int_{\Omega} |\nabla u_n|^2 \, dx} 
\\[2ex]
\leq & \multicolumn{2}{l}{\ds\lim_{n \to \infty} \Phi_{\varepsilon,\delta_n}^{\lambda_n} (u_n,u_{\Gamma,n}) - \varliminf_{n \to \infty} \int_{\Omega} f_{\delta_n}(\nabla u_n) \, dx - \varliminf_{n \to \infty} \int_{\Gamma} \big| \nabla_{\Gamma} (\varepsilon u_{\Gamma,n}) \big|^2 d\Gamma}
\\[2ex]
& & \ds- \varliminf_{ n \to \infty} \left( \int_{\Omega} B^{\lambda_n} (u_n) \, dx + \int_{\Gamma} B_{\Gamma}^{\lambda_n} (u_{\Gamma,n}) \, d\Gamma \right)
\\[2ex]
\leq & \multicolumn{2}{l}{\ds\Phi_{\varepsilon}(u,u_{\Gamma}) - \int_{\Omega} |{}\nabla u{}| \, dx - \int_{\Gamma} \bigl| \nabla_{\Gamma} (\varepsilon u_{\Gamma}) \bigr|^2 \, d\Gamma}
\\[2ex]
& & \ds - \left(\int_{\Omega} B (u) \, dx + \int_{\Gamma} B_{\Gamma} (u_{\Gamma}) \, d\Gamma \right) = \frac{\kappa^2}{2} \int_{\Omega} |\nabla u|^2 \,  dx.
\end{array}
\end{equation}
Having in mind \eqref{key1.10}, \eqref{key1.12-1} and the above calculation and the uniform convexity of $ L^2 $-based topologies, it is deduced that:
\begin{equation}\label{key1.13}
\left\{\begin{array}{l}
u_n \to u \mbox{ in $ H^1(\Omega) $, } 
\\[0.5ex]
\nabla u_n \to \nabla u \mbox{ in $ L^2(\Omega)^N $, } 
\end{array} \right.
\mbox{ as $ n \to \infty $.}
\end{equation}

In the meantime, from (A4) and \eqref{key1.13},
\begin{equation*}
\left|{}\nabla f_{\delta_n}(\nabla u_{n}){}\right|_{L^2(\Omega)^N}^2 \leq 2C_0 \left(\sup_{n \in \N} |{}\nabla u_n{}|_{L^2(\Omega)^N}^2 + \mathcal{L}^N(\Omega) \right) \mbox{, for any $ n \in \N $, }
\end{equation*}
which enables us to say
\begin{equation}\label{key1.14}
\nabla f_{\delta_n}(\nabla u_n) \to \nu_u \mbox{ weakly in $ L^2(\Omega)^N $ as $ n \to \infty $, for some $ \nu_u \in L^2(\Omega)^N $, }
\end{equation}
by taking a subsequence if necessary. 

In view of \eqref{key1.13}--\eqref{key1.14}, (Fact\,8), Remark \ref{exConvex}, and Remark \ref{Rem.MG} (Fact\,5), one can see that:
\begin{equation}\label{key1.15}
\nu_u \in \Sgn (\nabla u) \mbox{ a.e. in $ \Omega $.}
\end{equation}
Hence, letting $ n \to \infty $ in \eqref{key1.05} yields that:
\begin{equation}\label{key1.16}
\begin{array}{ll}
& \ds \int_{\Omega} \left( \nu_u + \kappa^2 \nabla u \right) \cdot \nabla z \, dx + \int_{\Gamma} \nabla_{\Gamma} (\varepsilon u_{\Gamma}) \cdot \nabla_{\Gamma} (\varepsilon z_{\Gamma}) \, d\Gamma + \int_{\Omega} \xi z \, dx + \int_{\Gamma} \xi_{\Gamma} z_{\Gamma} \, d\Gamma
\\[2ex]
= & \ds\int_{\Omega} ( w - u ) z \, dx + \int_{\Gamma} ( w_{\Gamma} - u_{\Gamma} ) z_{\Gamma} \, d\Gamma \mbox{, for any $ [z,z_{\Gamma}] \in \mathscr{V}_{\varepsilon} $.}
\end{array}
\end{equation}
In particular, taking any $ \varphi_0 \in H_0^1(\Omega) $ and putting $ [z,z_{\Gamma}] = [\varphi_0,0] $ in $ \mathscr{V}_{\varepsilon} $, 
\begin{equation*}
( w - u - \xi, \varphi_0 )_{L^2(\Omega)} = \int_{\Omega} (\nu_u + \kappa^2 \nabla u) \cdot \nabla \varphi_0 \, dx \mbox{, for any $ \varphi_0 \in H_0^1(\Omega) $.}
\end{equation*}
which implies:
\begin{equation}\label{key1.17}
-\mathrm{div} (\nu_u + \kappa^2 \nabla u) = w - u - \xi \in L^2(\Omega) \mbox{ in $ \mathscr{D}'(\Omega) $.}
\end{equation}
Furthermore, with Remark \ref{Rem.sOps} (Fact\,1)--(Fact\,3), \eqref{key1.16}--\eqref{key1.17} in mind, we can see that:
\begin{equation*}
\begin{array}{lll}
& \multicolumn{2}{l}{\ds(w_{\Gamma} - u_{\Gamma} - \xi_{\Gamma}, z_{\Gamma})_{L^2(\Gamma)}} 
\\[2ex]
= & \multicolumn{2}{l}{\ds\int_{\Omega} \left( \nu_u + \kappa^2 \nabla u \right) \cdot \nabla z \, dx - (w - u - \xi, z)_{L^2(\Omega)} + \int_{\Gamma} \nabla_{\Gamma} (\varepsilon u_{\Gamma}) \cdot \nabla_{\Gamma} (\varepsilon z_{\Gamma}) \, d\Gamma}
\\[2ex]
= & \multicolumn{2}{l}{{}_{H^{- \frac{1}{2}}(\Gamma)} \left\langle \left[\big( \nu_u + \kappa^2 \nabla u \big) \cdot n_{\Gamma} \right]_{_{\Gamma}}, z_{\Gamma} \right\rangle_{H^{\frac{1}{2}}(\Gamma)} + {}_{H^{- 1}(\Gamma)} \left\langle - {\mit \Delta}_{\Gamma} (\varepsilon u_{\Gamma}), \varepsilon z_{\Gamma} \right\rangle_{H^1(\Gamma)},}
\\[2ex]
\multicolumn{3}{c}{\mbox{for any $ [z,z_{\Gamma}] \in \mathscr{V}_{\varepsilon} $.}}
\end{array}
\end{equation*}
This identity leads to:
\begin{equation}\label{key1.18}
- {\mit \Delta}_{\Gamma} (\varepsilon^2 u_{\Gamma}) + \left[ \left( \nu_u + \kappa^2 \nabla u \right) \cdot n_{\Gamma} \right]_{_{\Gamma}} = w_{\Gamma} - u_{\Gamma} - \xi_{\Gamma} \in L^2(\Gamma) \mbox{ in $ H^{-1}(\Gamma) $. }
\end{equation}
As a consequence of \eqref{key1.10-1}, \eqref{key1.15}, \eqref{key1.17}--\eqref{key1.18}, we obtain {\it Claim $ \# $2}.

Now, by using {\it Claims $ \# $1--$ \# $2} and the maximality of $ \partial \Phi_{\varepsilon} $ in $ \mathscr{H} \times \mathscr{H} $, we can show the coincidence $ \partial \Phi_{\varepsilon} = \mathscr{A}_{\varepsilon} $ in $ \mathscr{H} \times \mathscr{H} $, and we conclude this Key-Lemma \ref{keyLem1}.
\hfill$ \Box $ 

\section{Proofs of Main Theorems}
\ \ \vspace{-3ex}

In this section, we will prove two Main Theorems by using the results of the previous sections. 
\\
\textbf{Proof of Main Theorem \ref{mainTh1}.}\ \ First, we show the item (I-1). In the Cauchy problem (CP)$_{\varepsilon}$, we see from (A3) and (Fact\,7) that:
\begin{equation*}
\Theta = [\theta,\theta_{\Gamma}] \in L^2(0,T;\mathscr{H}) \mbox{ and } U_0 = [u_0,u_{\Gamma,0}] \in \overline{D(\Phi_{\varepsilon})}.
\end{equation*}
Hence, by applying the general theories of evolution equations, e.g., \cite[\piegia{Theorem~4.1, p.~158}]{Barbu}, \cite[Theorem~3.6 \piegia{and Proposition~3.2}]{Brezis}, \cite[Section~2]{I-Y-K} and \cite[Theorem~1.1.2]{Kenmochi81}, we immediately have the existence and uniqueness of solution $ U = [u,u_{\Gamma}] \in L^2(0,T;\mathscr{H}) $ to (CP)$_{\varepsilon}$, such that:
\begin{equation*}
U \in C([0,T];\mathscr{H}) \piegia{{}\cap L^2(0,T; \mathscr{V}_{\varepsilon} )}\cap W_{\mathrm{loc}}^{1,2}((0,T];\mathscr{H}) \mbox{ and } \Phi_{\varepsilon} (U) \in L^1(0,T) \cap L_{\mathrm{loc}}^{\infty}((0,T]).
\end{equation*}
Also, there exists a positive constant $ C_1 $, independent of $ U_0 $ and $ \Theta $, such that:
\begin{equation}\label{MT1.1}
\begin{array}{c}
 \ds|U|_{C([0,T];\mathscr{H})\piegia{{}\cap L^2(0,T; \mathscr{V}_{\varepsilon} )}}^2 + \int_0^T \Phi_{\varepsilon} (U(t)) \, dt + \bigl| \sqrt{t} U' \bigr|_{L^2(0,T;\mathscr{H})}^{2} + \sup_{t \in (0,T)} t \Phi_{\varepsilon}(U(t))
\\[2ex]
\leq C_1 \left( 1 + |U_0|_{\mathscr{H}}^2 + |\Theta|_{L^2(0,T;\mathscr{H})}^2 \right). 
\end{array}
\end{equation}
Moreover, if $ U_0 \in D(\Phi_{\varepsilon})$, there exists a positive constant $ C_2 $, independent of $ U_0 $ and $ \Theta $, such that:
\begin{equation}\label{MT1.2}
\hspace{-2ex}|U'|_{L^2(0,T;\mathscr{H})}^2 + \sup_{t \in (0,T)} \Phi_{\varepsilon}(U(t)) \leq C_2 \left(1 + |U_0|_{\mathscr{H}}^2 + |\Theta|_{L^2(0,T;\mathscr{H})}^2 + \Phi_{\varepsilon}(U_0) \right).
\end{equation}

Now, Key-Lemma \ref{keyLem1} guarantees that the solution $ U = [u,u_{\Gamma}] $ to (CP)$_{\varepsilon}$ coincides with that to the system (ACE)$_{\varepsilon}$. Besides, in the light of  \eqref{Phi_nu} and (A1), the inequalities \eqref{est.1} and \eqref{est.2} directly follows from \eqref{MT1.1} and \eqref{MT1.2}, respectively.

Next, we show the item (I-2). For $ k= 1,2 $, let $ U^k := [u^k,u_{\Gamma}^k] $ be the two solutions to (CP)$_{\varepsilon}$ corresponds to forcing term $ \Theta^k := [\theta^k,\theta_{\Gamma}^k] \in L^2(0,T;\mathscr{H}) $ and initial term $ U_0^k = [u_0^k,u_{\Gamma,0}^k] \in D(\Phi_{\varepsilon}) $, respectively. Then, we can obtain the inequality \eqref{est.3} by using standard method: more precisely, by taking the difference between the two evolution equations, multiplying \piegia{it} by $ U^1(t) - U^2(t) $,  and applying Gronwall's lemma. \hfill$ \Box $
\medskip

\noindent
\textbf{Proof of Main Theorem \ref{mainTh2}.} The Main Theorem \ref{mainTh2} is proved by referring to the demonstration technique as in \cite[Theorem 2.7.1]{Kenmochi81}.

Let us set $ J_0 := [\varepsilon_0 - 1,\varepsilon_0 + 1] \cap [0, \infty) $. On this basis, we divide the proof in the following two steps.
\medskip

\noindent
{\it(Step\,1) The case when $ \{ \Phi_{\varepsilon} (\piegia{U_0^\varepsilon}) \}_{\varepsilon \in J_0} $ is bounded.} 

If $ \{ \Phi_{\varepsilon} (U_0^{\varepsilon}) \}_{\varepsilon \in J_0} $ is bounded, then the estimate \eqref{est.2} imply the following facts:
\begin{equation*}
\left\{\begin{array}{ll}
\ds\{ U^{\varepsilon} \}_{\varepsilon \in J_0} & \hspace{-2ex} \mbox{ is bounded in $ W^{1,2}(0,T;\mathscr{H}) \cap L^{\infty}(0,T;\mathscr{V}_0) $,}
\\[0.5ex] 
\ds\{ \varepsilon u_{\Gamma}^{\varepsilon} \}_{\varepsilon \in J_0} & \hspace{-2ex} \mbox{ is bounded in $ W^{1,2}(0,T;L^2(\Gamma)) \cap L^{\infty}(0,T;H^1(\Gamma)) $.}
\end{array} \right.
\end{equation*}
Therefore, applying general theories of compactness, such as Ascoli's theorem, we find a sequence $ \{ \varepsilon_n \}_{n=1}^{\infty} \subset J_0 $ and a limit point $ U = [u,u_{\Gamma}] \in W^{1,2}(0,T;\mathscr{H}) \cap L^{\infty}(0,T;\mathscr{V}_0) $, such that:
\begin{equation}\label{MT2.1}
\begin{array}{lll}
U^{\varepsilon_n} \to U & \multicolumn{2}{l}{\mbox{ in $ C([0,T];\mathscr{H}) $, weakly in $ W^{1,2}(0,T;\mathscr{H}) $}} 
\\[0.5ex]
& & \mbox{ and weakly-$*$ in $ L^{\infty}(0,T;\mathscr{V}_0) $ as $ n \to \infty $,}
\end{array}
\end{equation}
and
\begin{equation}\label{MT2.2}
\left\{\begin{array}{l}
\varepsilon_0 u_{\Gamma} \in W^{1,2}(0,T;L^2(\Gamma)) \cap L^{\infty}(0,T;H^1(\Gamma)),
\\[0.5ex]
\begin{array}{l} \hspace{-1ex} \varepsilon_n u_{\Gamma}^{\varepsilon_n} \to \varepsilon_0 u_{\Gamma} %& \multicolumn{2}{l}{\mbox{ in $ C([0,T];L^2(\Gamma)) $, weakly in $ W^{1,2}(0,T;L^2(\Gamma)) $}} 
%\\[0.5ex]
\mbox{ \ weakly-$*$ in $ L^{\infty}(0,T;H^1(\Gamma)) $, as $ n \to \infty $.}
\end{array}
\end{array} \right.
\end{equation}
Also, by Corollary \ref{cor.1} and Lemma \ref{auxLem} \piegia{with $S=(0,T)$, we have that}
\begin{equation*}
\piegia{\hat{\Phi}_{\varepsilon_n}} \to \piegia{\hat{\Phi}_{\varepsilon_0}} \mbox{ on $ L^2(0,T;\mathscr{H}) $, in the sense of \piegia{Mosco}, as $ n \to \infty $.}
\end{equation*}
From \eqref{MT2.1}, (A2), Remark \ref{Rem.Time-dep.} and Remark \ref{Rem.MG} (Fact\,5), it is seen that
\begin{equation}\label{MT2.3}
[-U' - \mathcal{G}(U) + \Theta^{\varepsilon_0}, U] \in \partial \piegia{\hat{\Phi}_{\varepsilon_0}} \mbox{ in $ L^2(0,T;\mathscr{H}) \times L^2(0,T;\mathscr{H}) $,}
\end{equation}
and
\begin{equation}\label{MT2.4}
\piegia{\hat{\Phi}_{\varepsilon_n}}(U^{\varepsilon_n}) \to \piegia{\hat{\Phi}_{\varepsilon_0}}(U) \mbox{ as $ n \to \infty $.}
\end{equation}
\piegia{Note that} \eqref{MT2.0-0}, \eqref{MT2.1}--\eqref{MT2.3} and Remark \ref{Rem.Time-dep.} (Fact\,4) \piegia{enable} us to say that $ U = [u,u_{\Gamma}] $ is a solution to the Cauchy problem (CP)$_{\varepsilon_0}$. So, due to the uniqueness of solutions, it must hold that:
\begin{equation}\label{MT2.5}
U = [u,u_{\Gamma}] = U^{\varepsilon_0} = [u^{\varepsilon_0},u_{\Gamma}^{\varepsilon_0}] \mbox{ in $ L^2(0,T;\mathscr{H}) $.}
\end{equation}

Furthermore, since (A1), \eqref{MT2.1}--\eqref{MT2.2} and \eqref{MT2.5} imply:
\begin{equation}\label{MT2.6}
\left\{\begin{array}{l}
\ds\varliminf_{n \to \infty} \frac{\kappa^2}{2} \int_0^T \! \int_{\Omega} |\nabla u^{\varepsilon_n}|^2 \, dx dt \geq \ds\frac{\kappa^2}{2} \int_0^T \! \int_{\Omega} |\nabla u^{\varepsilon_0}|^2 \, dx dt,
\\[2ex]
\ds\varliminf_{n \to \infty} \frac{1}{2} \int_0^T \! \int_{\Gamma} |\nabla_{\Gamma} (\varepsilon u_{\Gamma}^{\varepsilon_n})|^2 \, d\Gamma dt \geq \ds\frac{1}{2} \int_0^T \! \int_{\Gamma} |\nabla_{\Gamma} (\varepsilon u_{\Gamma}^{\varepsilon_0})|^2 \, d\Gamma  dt,
\\[2ex]
\ds \varliminf_{n \to \infty} \int_0^T \left( \int_{\Omega} \bigl( |\nabla u^{\varepsilon_n}| +B(u^{\varepsilon_n}) \bigr) \, dx + \int_{\Gamma} B_{\Gamma}(u_{\Gamma}^{\varepsilon_n}) \, d\Gamma \right) \, dt
\\[2ex]
\qquad\geq \ds\int_0^T \left( \int_{\Omega} \bigl( |\nabla u^{\varepsilon_0}| +B(u^{\varepsilon_0}) \bigr) \, dx + \int_{\Gamma} B_{\Gamma}(u_{\Gamma}^{\varepsilon_0}) \, d\Gamma \right)\, dt,
\end{array}\right.
\end{equation}
one can see that:
\begin{equation}\label{MT2.7}
\left\{\begin{array}{l}
\ds\lim_{n \to \infty} \int_0^T \! \int_{\Omega} | \nabla u^{\varepsilon_n} |^2 \, dx dt = \int_0^T \! \int_{\Omega} | \nabla u^{\varepsilon_0} |^2 \, dx dt,
\\[2ex]
\ds\lim_{n \to \infty} \int_0^T \! \int_{\Gamma} | \nabla_{\Gamma}(\varepsilon_n u_{\Gamma}^{\varepsilon_n}) |^2 \, d\Gamma dt = \int_0^T \! \int_{\Gamma} | \nabla_{\Gamma} (\varepsilon_0 u_{\Gamma}^{\varepsilon_0}) |^2 \, d\Gamma dt,
\end{array}\right.
\end{equation}
by applying (Fact\,0) with \eqref{MT2.4}--\eqref{MT2.6} in mind.

Now, taking into account \eqref{MT2.1}--\eqref{MT2.2}, \eqref{MT2.5} and \eqref{MT2.7}, and applying the uniform convexity of $ L^2 $-based topologies and the continuity of the trace operators, we obtain that:
\begin{equation}\label{MT2.7-1}
\left\{\begin{array}{lll}
u^{\varepsilon_n} \to u^{\varepsilon_0} \mbox{ in $ L^2(0,T;H^1(\Omega)) $, }\
\\[0.5ex]
\varepsilon_n u_{\Gamma}^{\varepsilon_n} \to \varepsilon_0 u_{\Gamma}^{\varepsilon_0} \mbox{ in $ L^2(0,T;H^1(\Gamma))$,}
\\[0.5ex]
u_{\Gamma}^{\varepsilon_n} = u^{\varepsilon_n}_{|_{\Gamma}} \to u_{\Gamma}^{\varepsilon_0} = u^{\varepsilon_0}_{|_{\Gamma}} \mbox{ in $ L^2(0,T;H^{\frac{1}{2}}(\Gamma)) $,}
\end{array}\right.
\mbox{ as $ n \to \infty $.}
\end{equation}
\piegia{Therefore,} \eqref{MT2.1}, \eqref{MT2.5} and \eqref{MT2.7-1} are sufficient to verify \eqref{MT2.0-1}--\eqref{MT2.0-2}.
\medskip

\noindent
{\it(Step\,2) The case when $ \{ \Phi_{\varepsilon}(U_0^{\varepsilon}) \}_{\varepsilon \in J_0} $ is unbounded.}

Let $ \rho \in (0,1) $ be an arbitrary constant. Then, by the assumption for $ \{ U_0^{\varepsilon} \}_{\varepsilon \geq 0} $, \piegia{for any sequence $\{ \varepsilon_n\}_{n=1}^\infty \subset J_0$ converging to $\varepsilon_0$} we find a large number $ n_1(\rho) \in \N $, such that:
\begin{equation}\label{MT2.10}
\left| U_0^{\varepsilon_n} - U_0^{\varepsilon_0} \right|_{\mathscr{H}} \leq \rho \mbox{, for any $ n \geq n_1(\rho) $ }
\end{equation}
Also, since $ U_0^{\varepsilon_0} \in \overline{D(\Phi_{\varepsilon_0})} $, we find a function $ W_{0,\rho} \in D(\Phi_{\varepsilon_0}) $, such that:
\begin{equation}\label{MT2.11}
\left| U_0^{\varepsilon_0} - W_{0,\rho} \right|_{\mathscr{H}} \leq \rho.
\end{equation}
Additionally, by Corollary \ref{cor.1}, there exists a sequence $ \{ W_{0,\rho}^n \in D(\Phi_{\varepsilon_n}) \}_{n=1}^{\infty} \subset \mathscr{H} $, such that:
\begin{equation}\label{MT2.11-1}
\left\{ \begin{array}{l}
W_{0,\rho}^n \to W_{0,\rho} \mbox{ in $ \mathscr{H} $,}
\\[0.5ex]
\Phi_{\varepsilon_n}(W_{0,\rho}^n) \to \Phi_{\varepsilon_0}(W_{0,\rho}),
\end{array}\right.
\mbox{as $ n \to \infty $,}
\end{equation}
and in particular, there exists a large number $ n_2(\rho) \in \N $, with $ n_2(\rho) \geq n_1(\rho) $, such that:
\begin{equation}\label{MT2.12}
| W_{0, \rho}^n - W_{0,\rho} |_{\mathscr{H}} \leq \rho \mbox{, for any $ n \geq n_2({\rho}) $.}
\end{equation}
From \eqref{MT2.10}--\eqref{MT2.11} and \eqref{MT2.12}, it follows that:
\begin{equation}\label{MT2.13}
\begin{array}{rcl}
| U_0^{\varepsilon_n} - W_{0,\rho}^n |_{\mathscr{H}} & \leq & | U_0^{\varepsilon_n} - U_0^{\varepsilon_0} |_{\mathscr{H}} + | U_0^{\varepsilon_0} - W_{0, \rho} |_{\mathscr{H}} + | W_{0,\rho} - W_{0,\rho}^n |_{\mathscr{H}}
\\[2ex]
& \leq & 3 \rho \mbox{, for any $ n \geq n_2(\rho) $.}
\end{array}
\end{equation}
Based on \piegia{this}, let $ W_{\rho} \in W^{1,2}(0,T;\mathscr{H}) \cap L^{\infty}(0,T;\mathscr{V}_{\varepsilon_0}) $ be the solution to (CP)$_{\varepsilon_0}$, corresponding to the forcing term $ \Theta^{\varepsilon_0} \in L^2(0,T;\mathscr{H}) $ and the initial data $ W_{0,\rho} \in D(\Phi_{\varepsilon_0}) $. As well as, for any $ n \in \N $, let $ W_{\rho}^n \in W^{1,2}(0,T;\mathscr{H}) \cap L^{\infty}(0,T;\mathscr{V}_{\varepsilon_n}) $ be the solution to (CP)$_{\varepsilon_n}$, corresponding to the forcing term $ \Theta^{\varepsilon_n} \in L^2(0,T;\mathscr{H}) $ and the initial data $ W_{0,\rho}^n \in D(\Phi_{\varepsilon_n}) $. Then, by applying the result of the previous {(Step\,1)}, we have:
\begin{equation}\label{MT2.14}
W_{\rho}^n \to W_{\rho} \mbox{ in $ C([0,T];\mathscr{H}) $ as $ n \to \infty $.}
\end{equation}
Besides, from \eqref{est.3} and \eqref{MT2.13}--\eqref{MT2.14}, one can see that:
\begin{equation*}
\begin{array}{lll}
& \multicolumn{2}{l}{\ds\varlimsup_{n \to \infty} \left| U^{\varepsilon_n} - U^{\varepsilon_0} \right|_{C([0,T];\mathscr{H})}}
\\[2ex]
\leq & \multicolumn{2}{l}{\ds\varlimsup_{n \to \infty} \left( \left| U^{\varepsilon_n} - W_{\rho}^n \right|_{C([0,T];\mathscr{H})} + \left| W_{\rho}^n - W_{\rho} \right|_{C([0,T];\mathscr{H})} + \left| W_{\rho} - U^{\varepsilon_0} \right|_{C([0,T];\mathscr{H})} \right)}
\\[2ex]
\leq & \multicolumn{2}{l}{\ds\sqrt{C_3} \left( \varlimsup_{n \to \infty} \left| U_0^{\varepsilon_n} - W_{0,\rho}^n \right|_{\mathscr{H}} + \left| W_{0,\rho} - U_0^{\varepsilon_0} \right|_{\mathscr{H}} \right) + \lim_{n \to \infty} \left| W_{\rho}^n - W_{\rho} \right|_{C([0,T];\mathscr{H})}}
\\[2ex]
\leq & \multicolumn{2}{l}{4 \sqrt{C_3} \rho.}
\end{array}
\end{equation*}
Since $ \rho \in (0,1) $ is arbitrary, the above inequality implies that:
\begin{equation}\label{MT2.19}
U^{\varepsilon_n} \to U^{\varepsilon_0} \mbox{ in $ C([0,T];\mathscr{H}) $ as $ n \to \infty $.}
\end{equation}

Now, our remaining task will be to verify the convergences \eqref{MT2.0-1}--\eqref{MT2.0-2} under the unbounded situation of $ \{ \Phi_{\varepsilon_n}(U_0^{\varepsilon_n}) \}_{n=1}^{\infty} $. To this end, we first invoke \eqref{est.2} and \eqref{MT2.11-1} to check the existence of a constant $ K(\rho) $, depending on $ \rho \in (0,1) $, such that:
\begin{equation}\label{MT2.20}
\begin{array}{c}
\ds| (W_{\rho})' |_{L^2(0,T;\mathscr{H})}^2 + | (W_{\rho}^n)' |_{L^2(0,T;\mathscr{H})}^2 + \sup_{t \in [0,T]} \Phi_{\varepsilon_0}(W_{\rho}(t)) + \sup_{t \in [0,T]} \Phi_{\varepsilon_n}(W_{\rho}^n(t)) 
\\[2ex]
\leq K(\rho) \mbox{, for any $ n \in \N $.}
\end{array}
\end{equation}
On this basis, let us consider a sequence $ \{ \mathcal{P}_n \}_{n=1}^{\infty} \subset C([0,T]) $ of functions, given as:
\begin{equation*}
\varsigma \in [0,T] \mapsto \mathcal{P}_n(\varsigma) := \int_0^{\varsigma} \Phi_{\varepsilon_n}(W_{\rho}^n(t)) \, dt \in [0,\infty) \mbox{, for any $ n \in \N $.}
\end{equation*}
Then, by applying a similar method to show \eqref{MT2.4}, we have:
\begin{equation}\label{MT2.21}
\mathcal{P}_n(\varsigma) \to \mathcal{P}(\varsigma) := \int_0^{\varsigma} \Phi_{\varepsilon_0}(W_{\rho}(t)) \, dt \mbox{ as $ n \to \infty $, for any $ \varsigma \in [0,T] $.}
\end{equation}
Also, from \eqref{MT2.20}, it is seen that:
\begin{equation}\label{MT2.22}
\left\{ \frac{d}{d\varsigma} \mathcal{P}_n \right\}_{n=1}^{\infty} = \left\{ \Phi_{\varepsilon_n} (W_{\rho}^n) \right\}_{n=1}^{\infty} \mbox{ is bounded in $ L^{\infty}(0,T) $.}
\end{equation}
By \eqref{MT2.21}--\eqref{MT2.22} and Ascoli's theorem, we may suppose that:
\begin{equation*}
\mathcal{P}_n \to \mathcal{P} \mbox{ in $ C([0,T]) $ as $ n \to \infty $,}
\end{equation*}
by taking a subsequence if necessary, and more precisely, we find a large number $ n_*(\rho) $, independent of $ \varsigma \in [0,T] $, such that:
\begin{equation}\label{MT2.23}
\begin{array}{c}
\ds n_*(\rho) \geq n_2(\rho) \mbox{ and } \left| \int_0^{\varsigma} \Phi_{\varepsilon_n}(W_{\rho}^n(t)) \, dt - \int_0^{\varsigma} \Phi_{\varepsilon_0}(W_{\rho}(t)) \, dt \right| < \rho,
\\[2ex]
\mbox{for any $ n \geq n_*(\rho) $ and any $ \varsigma \in [0,T] $.}
\end{array}
\end{equation}

In the meantime, for the sequence of solutions $ \{ U^{\varepsilon_n} \}_{n=1}^{\infty} $, it is easily seen that:
\begin{equation*}
\begin{array}{lll}
\multicolumn{3}{l}{\ds\left( (U^{\varepsilon_n} - W_{\rho}^n)' , U^{\varepsilon_n} - Z \right)_{L^2(0,\varsigma;\mathscr{H})} + \int_0^{\varsigma} \Phi_{\varepsilon_n}(U^{\varepsilon_n}(t)) \, dt}
\\[2ex] 
& \multicolumn{2}{l}{\ds\leq \int_0^{\varsigma} \Phi_{\varepsilon_n}(Z(t)) \, dt + \left( \Theta^{\varepsilon_n} - \mathcal{G}(U^{\varepsilon_n}) - (W_{\rho}^n )', U^{\varepsilon_n} - Z \right)_{L^2(0,\varsigma;\mathscr{H})},}
\\[2ex]
\multicolumn{3}{c}{\mbox{for any $ Z \in L^2(0,\varsigma;\mathscr{V}_{\varepsilon_n}) $ and any $ n \in \N $.}}
\end{array}
\end{equation*}
So, putting $ Z = W_{\rho}^n $, and using \eqref{est.3}--\eqref{MT2.0-0}, \eqref{MT2.13}, \eqref{MT2.19}--\eqref{MT2.20}, \eqref{MT2.23} and (A2), we infer that:
\begin{equation}\label{MT2.24}
\begin{array}{llll}
& \multicolumn{3}{l}{\ds\int_0^{\varsigma} \Phi_{\varepsilon_n}(U^{\varepsilon_n}(t)) \, dt \leq \frac{1}{2} | U_0^{\varepsilon_n} - W_{0,\rho}^n |_{\mathscr{H}}^2 + \int_0^{\varsigma} \Phi_{\varepsilon_n}(W_{\rho}^n(t)) \, dt} 
\\[2ex]
& & & + \ds | \Theta^{\varepsilon_n} - \mathcal{G}(U^{\varepsilon_n}) - ( W_{\rho}^n )' |_{L^2(0,\varsigma;\mathscr{H})} | U^{\varepsilon_n} - W_{\rho}^n |_{L^2(0,\varsigma;\mathscr{H})}
\\[2ex]
\leq & \multicolumn{3}{l}{\ds\frac{9}{2} \rho^2 + \left( \int_0^{\varsigma} \Phi_{\varepsilon_0}(W_{\rho}(t)) \, dt + \rho \right)} 
\\[2ex]
& & \multicolumn{2}{l}{+ \ds\sqrt{C_3 \varsigma} | U_0^{\varepsilon_n} - W_{0,\rho}^n |_{\mathscr{H}}\left (|\Theta^{\varepsilon_n} - \mathcal{G}(U^{\varepsilon_n})|_{L^2(0,\varsigma;\mathscr{H})} + | \bigl( W_{\rho}^n \bigr)' |_{L^2(0,\varsigma;\mathscr{H})} \right)}
\\[2ex]
\leq & \multicolumn{3}{l}{\ds\frac{11}{2} \rho + \int_0^{\varsigma} \Phi_{\varepsilon_0} (W_{0,\rho}(t)) \, dt}
\\[2ex]
& & \multicolumn{2}{l}{+ 3\rho \sqrt{C_3 \varsigma} \left(\sup_{n \in \N} |\Theta^{\varepsilon_n} - \mathcal{G}(U^{\varepsilon_n})|_{L^2(0,\varsigma;\mathscr{H})} + \sqrt{K(\rho)} \right),}
\\[2ex]
\multicolumn{4}{c}{\mbox{for any $ \varsigma \in [0,T] $ and any $ n \in \N $.}}
\end{array}
\end{equation}
Here, let us take a constant $ \varsigma_*(\rho) \in (0,T) $, so small to satisfy that:
\begin{equation}\label{MT2.25}
\begin{array}{l}
\ds\int_0^{\varsigma_*(\rho)} \Phi_{\varepsilon_0}(W_{0,\rho}(t)) \, dt 
\\[2ex]
\qquad \ds+ 3\rho \sqrt{C_3 \varsigma_*(\rho)} 
\left( \sup_{n \in \N} | \Theta^{\varepsilon_n} - \mathcal{G}(U^{\varepsilon_n}) 
|_{\piegia{L^2(0,\varsigma_*(\rho);\mathcal{H})}} + K(\rho) \right) < \frac{\rho}{2}.
\end{array}
\end{equation}
Then, having in mind \eqref{MT2.19}, \eqref{MT2.24}--\eqref{MT2.25}, Corollary \ref{cor.1} and Fatou's lemma, it is observed that:
\begin{equation}\label{MT2.26}
\begin{array}{lll}
\ds\int_0^{\varsigma_*(\rho)} \Phi_{\varepsilon_0} (U^{\varepsilon_0}(t)) \, dt & \leq & \ds\varliminf_{n \to \infty} \int_0^{\varsigma_*(\rho)} \Phi_{\varepsilon_n} (U^{\varepsilon_n}(t)) \, dt
\\[2ex]
& \leq & \ds \sup_{n \in \N} \int_0^{\varsigma_*(\rho)} \Phi_{\varepsilon_n}(U^{\varepsilon_n}(t)) \, dt \leq 6 \rho.
\end{array}
\end{equation}

Finally, by \eqref{est.1}, one can see that the sequence $ \{ U^{\varepsilon_n} \}_{n=1}^{\infty} $ is bounded in $ W^{1,2}(\varsigma_*(\rho),T; $ $ \mathscr{H}) \cap L^{\infty}(\varsigma_*(\rho),T;\mathscr{V}_0) $. So, we can apply a similar arguments to obtain \eqref{MT2.4}, and we can show that:
\begin{equation}\label{MT2.27}
\int_{\varsigma_*(\rho)}^T \Phi_{\varepsilon_n}(U^{\varepsilon_n}(t)) \, dt \to \int_{\varsigma_*(\rho)}^T \Phi_{\varepsilon_0}(U^{\varepsilon_0}(t)) \, dt \mbox{ as $ n \to \infty $.}
\end{equation}

In view of \eqref{MT2.19}, \eqref{MT2.26}--\eqref{MT2.27}, we can say that:
\begin{equation}\label{MT2.28}
\left\{\begin{array}{l}
\mbox{$ U^{\varepsilon_0} \in L^2(0,T;\mathscr{V}_0)$, $ \{ U^{\varepsilon_n} \}_{n=1}^{\infty} $ is bounded in $ L^2(0,T;\mathscr{V}_0) $,}
\\[0.5ex]
\mbox{$ U^{\varepsilon_n} \to U^{\varepsilon_0} $ weakly in $ L^2(0,T;\mathscr{V}_0) $ as $ n \to \infty $,}
\end{array}\right.
\end{equation}
and
\begin{equation*}
\begin{array}{lll}
& \multicolumn{2}{l}{\ds\varlimsup_{n \to \infty} \left| \int_0^{T} \Phi_{\varepsilon_n}(U^{\varepsilon_n}(t)) \, dt - \int_0^{T} \Phi_{\varepsilon_0}(U^{\varepsilon_0}(t)) \, dt \right|}
\\[2ex]
\leq & \multicolumn{2}{l}{\ds\sup_{n \in \N} \int_0^{\varsigma_*(\rho)} \Phi_{\varepsilon_n}(U^{\varepsilon_n}(t)) \, dt + \int_0^{\varsigma_*(\rho)} \Phi_{\varepsilon_0}(U^{\varepsilon_0}(t)) \, dt}
\\[2ex]
& & + \ds\lim_{n \to \infty}\left|\int_{\varsigma_*(\rho)}^{T} \Phi_{\varepsilon_n}(U^{\varepsilon_n}(t)) \, dt - \int_{\varsigma_*(\rho)}^{T} \Phi_{\varepsilon_0}(U^{\varepsilon_0}(t)) \, dt \right|
\\[2ex]
\leq & \multicolumn{2}{l}{12 \rho.}
\end{array}
\end{equation*}
Since $ \rho \in (0,1) $ is arbitrary, the above inequality implies:
\begin{equation}\label{MT2.29}
\lim_{n \to \infty} \int_0^T \Phi_{\varepsilon_n}(U^{\varepsilon_n}(t)) \, dt = \int_{0}^T \Phi_{\varepsilon_0}(U^{\varepsilon_0}(t)) \, dt.
\end{equation}
By virtue of \eqref{MT2.28}--\eqref{MT2.29}, we can apply a similar method  to derive \eqref{MT2.7-1}, and we obtain that:
\begin{equation}\label{MT2.30}
\left\{\begin{array}{lll}
u^{\varepsilon_n} \to u^{\varepsilon_0} \mbox{ in $ L^2(0,T;H^1(\Omega)) $, }\
\\[2ex]
\varepsilon_n u_{\Gamma}^{\varepsilon_n} \to \varepsilon_0 u_{\Gamma}^{\varepsilon_0} \mbox{ in $ L^2(0,T;H^1(\Gamma))$,}
\\[2ex]
u_{\Gamma}^{\varepsilon_n} = u^{\varepsilon_n}_{|_{\Gamma}} \to u_{\Gamma}^{\varepsilon_0} = u^{\varepsilon_0}_{|_{\Gamma}} \mbox{ in $ L^2(0,T;H^{\frac{1}{2}}(\Gamma)) $,}
\end{array}\right.
\mbox{ as $ n \to \infty $.}
\end{equation}
\piegia{Hence,} \eqref{MT2.19} and \eqref{MT2.30} \piegia{imply} the conclusive convergences \eqref{MT2.0-1}--\eqref{MT2.0-2}.
\hfill$ \Box $

\end{document}